\documentclass[11pt, english]{amsart}
\usepackage{amssymb,amsthm,amsmath,amstext,enumitem}
\usepackage[utf8]{inputenc}
\usepackage[left=2.4cm,top=2.4cm,right=2.4cm,bottom=2.4cm]{geometry}
\usepackage[all]{xy}
\usepackage{eucal}
\usepackage[pdfusetitle]{hyperref}
\usepackage[only,llbracket,rrbracket,sslash]{stmaryrd}
\usepackage{mathtools}
\usepackage{dsfont}
\usepackage[usenames,dvipsnames]{xcolor}

\makeatletter
\renewcommand\paragraph{%
  \@startsection{paragraph}{4}{\z@}%
    {1.5ex \@plus .5ex \@minus .2ex}%
    {-0.5em}%
    {\normalfont\normalsize\itshape}%
}
\makeatother

\theoremstyle{plain}
\newtheorem{theorem}{Theorem}[section]
\newtheorem{proposition}[theorem]{Proposition}
\newtheorem{lemma}[theorem]{Lemma}

\newtheorem{corollary}[theorem]{Corollary}
\theoremstyle{definition}
\newtheorem{definition}[theorem]{Definition}
\newtheorem{remark}[theorem]{Remark}
\newtheorem{example}[theorem]{Example}

\DeclareMathOperator{\Spec}{Spec}
\DeclareMathOperator{\rk}{rk}

\DeclareMathOperator{\Hom}{Hom}
\DeclareMathOperator{\Der}{Der}
\DeclareMathOperator{\Coker}{coker}
\DeclareMathOperator{\depth}{depth}

\DeclareMathOperator{\id}{\mathds{1}}
\DeclareMathOperator{\im}{Im}
\DeclareMathOperator{\chiV}{\chi}

\def\bA{{\mathbf{A}}}
\def\bM{{\mathbf{M}}}
\def\bU{{\mathbf{U}}}
\def\Derrho{{\mathrm{Der_\bU^{\rho}}}}
\def\Derrhoad{{\mathrm{Der_\bU^{\mathrm{ad}}}}}
\def\bE{{\mathbf{E}}}
\def\bR{{\mathbf{R}}}
\def\bS{{\mathbf{S}}}
\def\bH{{\mathbf{H}}}
\def\bF{{\mathbf{F}}}
\def\bx{{\mathbf{x}}}

\renewcommand{\AA}{{\mathds{A}}}
\def\kk{{\mathds{k}}}
\def\ZZ{{\mathds{Z}}}

\def\PP{{\mathds{P}}}
\def\GG{{\mathds{G}}}
\def\VV{{\mathds{V}}}

\def\SL{\mathrm{SL}}
\def\GL{\mathrm{GL}}
\def\Spin{\mathrm{Spin}}
\def\LL{{L}}
\def\red{{\mathrm{red}}}
\def\Sp{\mathrm{Sp}}

\def\cT{{\mathcal{T}}}
\def\cO{{\mathcal{O}}}

\def\fg{{\mathfrak {g}}}
\def\fu{{\mathfrak {u}}}
\def\fso{{\mathfrak {so}}}
\def\gl{{\mathfrak {gl}}}
\def\fh{{\mathfrak {h}}}
\def\fe{{\mathfrak {e}}}
\def\fsl{{\mathfrak {sl}}}
\def\fsp{{\mathfrak {sp}}}
\def\ft{{\mathfrak {t}}}
\def\fL{{\mathfrak {L}}}

\def\Gmod{\operatorname{Rep}(G)}
\def\Gvmod{\operatorname{Rep}(G_v)}
\def\Wmod{\operatorname{Rep}(W)}

\def\Exp{{\mathrm{Exp}}}
\def\Coexp{{\mathrm{Coexp}}}

\begin{document}

\author{Vladimiro Benedetti}
\address{Vladimiro Benedetti, 
Université Côte d’Azur, CNRS – Laboratoire J.-A. Dieudonné, Parc Valrose, F-06108 Nice Cedex 2, France}
\email{vladimiro.benedetti@univ-cotedazur.fr}

\author{Daniele Faenzi}
\address{Daniele Faenzi, 
Université Bourgogne Europe, CNRS, IMB UMR 5584, F-21000 Dijon, France
}
\email{daniele.faenzi@u-bourgogne.fr}

\author{Simone Marchesi}
\address{Simone Marchesi, 
 Facultat de Matemàtiques i Informàtica, Universitat de Barcelona, Gran
Via de les Corts Catalanes 585, 08007 Barcelona, Spain}
\address{Centre de Recerca Matemàtica Edifici C, Campus Bellaterra, 08193 Bellaterra, Spain}
\email{marchesi@ub.edu}

\title{Derivations tangent to polar discriminants}

\begin{thanks}{V. B. and D. F. were partially supported by FanoHK ANR-20-CE40-0023, SupToPhAG/EIPHI ANR-17-EURE-0002, Bridges ANR-21-CE40-0017. This work was started during the stay of D. F. at Osaka University, funded by JSPS Invitational Fellowship for Research in Japan (short-term) ID S24043. S. M. was partially supported by PID2024-157142NB-I00 and by the Spanish State Research Agency, through the Severo Ochoa
and María de Maeztu Program for Centers and Units of Excellence in R\&D (CEX2020-001084-M).
}\end{thanks}



\hypersetup{
  pdftitle={Derivations tangent to polar discriminants},
  pdfauthor={Vladimiro Benedetti, Daniele Faenzi, Simone Marchesi},
  pdfkeywords={Polar representations, discriminants, logarithmic derivations, free divisors}
}

\begin{abstract}
We describe the module of logarithmic derivations along the
discriminant of a stable polar representation of a reductive
algebraic group. It splits into derivations along the fibres of the
quotient map and a free summand lifted from logarithmic derivations
along the associated branch divisor. This gives a uniform extension
of previous results on full adjoint discriminants. We prove that, for visible stable polar
actions of semisimple groups with finite generic stabilizer, the discriminant is a free divisor. For stable periodic
gradings of arbitrary order with semisimple degree-zero group, we
construct explicit Saito matrices from the infinitesimal action and
the gradients of basic invariants, which yields determinantal
equations for several exceptional discriminants. In particular, starting with exceptional Lie algebras, we provide examples
of free discriminants of degrees 60, 120, 126 and 240 for inner gradings and a degree-72 example for an outer one.
\end{abstract}

\maketitle

\section{Introduction}

We work over an algebraically closed field \(\kk\) of
characteristic zero.
A reduced hypersurface \(D=\VV(F)\) in an \(\ell\)-dimensional
\(\kk\)-vector space \(V\), with reduced equation
\(F\in\bU=\kk[V]\), is said to be \textit{free} if its module of
logarithmic derivations \(\Der_\bU(-\log F)\), consisting of derivations
\(\vartheta\) of \(\bU\) such that \(F\) divides \(\vartheta(F)\), is a free
\(\bU\)-module. Saito introduced this notion in
\cite{saito:logarithmic}, where he developed a celebrated criterion:
\(D\) is free if and only if there exist \(\ell\) derivations in
\(\Der_\bU(-\log F)\) whose matrix of coefficients (the \textit{Saito
matrix}) has determinant equal to a unit times \(F\).
Free divisors play a central role in the study of singularities and their unfoldings, see e.g. \cite{loojenga:semi-universal, buchweitz-et-al}.
In the framework of hyperplane arrangements, freeness is one of the most  studied properties, see e.g. \cite{abe:divisionally, yoshinaga:free, orlik-terao:arrangements, terao:exponents}. 
As it will be crucial for us, if \(D\) is the set of reflecting hyperplanes of a finite group \(W\) generated by pseudo-reflections, then \(D\) is free and the degrees of a free basis of logarithmic derivations are precisely the coexponents of \(W\), see \cite{orlik-solomon:unitary,terao:free-unitary,orlik-terao:arrangements}.

From a different perspective, interesting hypersurfaces arise in the context of invariant theory, the connection with arrangements being provided by Chevalley's theorem.
For a reductive group \(G\) with Lie algebra \(\fg\), Cartan subalgebra \(\fh\) and Weyl group
\(W\), it gives
\(
\fg\sslash G\simeq\fh/W
\)
and
\(
\kk[\fg]^G\simeq\kk[\fh]^W.
\)
The reduced preimage of the branch divisor under
\(\fh\to\fh/W\) is the reflection arrangement of \(W\), which is free.
The preimage of the branch divisor under the adjoint quotient
\(\fg\to\fg\sslash G\) is a \(G\)-invariant hypersurface in \(\fg\), which
we call the \emph{full adjoint discriminant}.
The study of adjoint discriminants has been recently extended to the case of graded Lie algebras in \cite{benedetti-manivel:discriminants}.
These discriminants are connected to moduli spaces of abelian varieties of small dimension and curves of low genus, according to recent work of many authors, see \cite{GSW, gruson-sam, sam-rains1, Thorne, sam-rains2,BBFM1, laga:genus3, laga-romano} to mention a few.

The goal of this paper is to study the module of logarithmic derivations
of invariant hypersurfaces of this kind, in a substantially more general framework, focusing on freeness and, more generally, on the projective dimension of \(\Der_\bU(-\log F)\).
The natural setting is that of
polar representations, introduced by Dadok and Kac in
\cite{dadok-kac}, see \S\ref{section:polar}.
They include Vinberg theta representations arising from periodic
gradings of Lie algebras, and in particular adjoint representations.
Our previous work
\cite{benedetti-faenzi-marchesi:adjoint-arxiv} treated adjoint
discriminants for simple groups of simply-laced type using
projective-duality methods coupled with a case-by-case analysis. The approach developed here is uniform
and applies in the broader setting considered below.

Following Dadok and Kac \cite{dadok-kac}, a representation \(V\) of a
reductive group \(G\) is called \emph{polar} if there exists a
semisimple element \(v\in V\) such that
\(
\dim(V\sslash G)=\dim_{\kk}(\fh_v),
\)
where 
\(
\fh_v=\{x\in V\mid\fg\cdot x\subseteq\fg\cdot v\}.
\)
The space \(\fh=\fh_v\) is called a \emph{Cartan subspace}; see
\S\ref{section:polar}. For such a representation, we have a Chevalley-type theorem, asserting that
\(
V\sslash G\simeq\fh/W
\) is an affine space,
where the Weyl group
\(W=N_G(\fh)/Z_G(\fh)
\)
is a finite pseudo-reflection group. Writing
\(\bU=\kk[V]\), \(\bS=\kk[\fh]\), and \(\bR=\bU^G\), one has
\(\bR\simeq\bS^W\). We set
\(\Coexp(V)=\Coexp(W)\); see \S\ref{sec:exponents}.
Let \(H_1,\ldots,H_N\) be the reflecting hyperplanes of \(W\), let
\(\delta_j\) be a linear equation of \(H_j\), and let \(m_j\) be the
order of the pointwise stabilizer of \(H_j\), for \(1\le j\le N\). Then
\(
\delta=\prod_{j=1}^N\delta_j^{m_j}\in\bS^W\simeq\bR.
\)
Viewed in \(\bS\), its divisor is \(\sum_{j=1}^Nm_jH_j\), whereas, viewed in
\(\bR\), it is a reduced equation of the branch divisor of
\(\fh\to\fh/W\). We denote by \(\Delta\in\bU^G\) the corresponding invariant
polynomial and set
\(
D=\VV(\Delta)\subset V.
\)
We call \(\Delta\) the \emph{polar discriminant equation} and \(D\) the
\emph{polar discriminant divisor}; see \S\ref{sec-polardiscriminant}.

The quotient morphism \(\chiV:V\to V\sslash G\) induces, by restriction
to \(\bR\subset\bU\),
the differential
\(
d\chiV:\Der_\bU\to\Der_{\bR,\bU}
\) defined by sending a derivation \(\vartheta\) to \(\vartheta|_\bR\).
The module
\(
\Der_\bU^{\chiV}
=\ker(d\chiV)
\)
consists of vector fields tangent to the fibres of \(\chiV\), see
\S\S\ref{section:variant} and \ref{sec-actionderivation}.
Our first main result decomposes \(\Der_\bU(-\log\Delta)\) into
the derivations along the fibres of \(\chiV\) and a free summand pulled
back from the branch divisor, whose reduced preimage in \(\fh\) is precisely given by the reflection arrangement associated with the Weyl group \(W\).
For a general \(v\in\fh\), let \(G_v\) be its stabilizer. One has
\(\fh\subseteq V^{G_v}\).

\begin{theorem} \label{main:general}
Let \(V\) be a stable polar representation of a reductive group \(G\). Then \(\Der_\bU^G\) is a free \(\bR\)-module of rank
\(\dim_{\kk}(V^{G_v})\), where \(v\in\fh\) is general, and  \(\Der_\bR(-\log\delta) \otimes_\bR \bU\) is a free \(\bU\)-module, isomorphic to \(\bigoplus_{c \in \Coexp(V)} \bU(-c)\).
Moreover, if the restriction map
\[
\pi_G:\Der_\bU^G \xrightarrow{\simeq} \Der_\bR(-\log\delta),
\qquad
\eta\mapsto\eta|_\bR.
\]
is surjective, then there is a direct-sum decomposition
\[
\Der_\bU(-\log \Delta) \simeq \Der_{\bU}^{{\chiV}} \oplus  \Der_\bR(-\log\delta) \otimes_\bR \bU.
\]
In particular, if
\(\fh = V^{G_v}\),  then $\pi_G$ is an isomorphism.
\end{theorem}

Stability here means that a general \(G\)-orbit is closed. For a polar representation this is equivalent to the fact that \(V\) decomposes as \(V=\fh\oplus\fg\cdot\fh\); see Remark~\ref{rmk-stability}.
The theorem has a sheaf-theoretic counterpart, more geometric in flavor, that we spell out in Corollary~\ref{main-corollary}.
Whenever $\pi_G$ is surjective, we will say that we have a \textit{surjective restriction of invariant derivations}.

Theorem~\ref{main:general} applies in particular to the full adjoint
discriminant of a semisimple algebraic group \(G\); see
Theorem~\ref{theorem:adjoint}. It gives a decomposition
\[
\Der_\bU(-\log\Delta)
\simeq
\Derrhoad\oplus
\bigoplus_{e\in\Exp(\fg)}\bU(-e),
\]
where \(\Derrhoad\) has projective dimension one and admits an explicit
presentation in terms of the gradients of the basic invariants of
\(G\). For simple groups of simply-laced type, this was established in
\cite{benedetti-faenzi-marchesi:adjoint-arxiv} by projective-duality
methods inspired by \cite{weyman:syzygies,GKZ:discriminants}; here it
follows from the uniform approach developed above.

A further consequence concerns polar actions with finite generic
stabilizer; see Theorem~\ref{thm_finite_stab}. Recall that \(V\) is
\emph{visible} if every fibre of \(\chiV\) contains finitely many
\(G\)-orbits. The Vinberg theta representations considered below are
visible.

\begin{theorem}
Let \(V\) be a visible stable polar representation of a semisimple group \(G\) with finite generic stabilizer and surjective restriction of invariant derivations. Then the polar discriminant \(\Delta\) is a free divisor, more precisely:
\[
\Der_\bU(-\log \Delta)\simeq\fg\otimes \bU(-1) \oplus \bigoplus_{c \in \Coexp(V)} \bU(-c).
\]
\end{theorem}

As main application, Saito's criterion gives explicit determinantal
representations of discriminants for stable periodic gradings of arbitrary order.
Let
\(\LL\) be a connected semisimple algebraic group with Lie algebra
\(\fL\), and assume that
\(
\fL=\bigoplus_{a\in\ZZ/m\ZZ}\fL_a
\)
is the grading defined by an automorphism \(\theta\) of \(\LL\) of
order \(m\). Set
\(
G=(\LL^\theta)^\circ
\)
and
\(
V=\fL_1.
\)
Let
\(
\Phi_1,\ldots,\Phi_r\in\kk[\fL]^\LL
\)
be homogeneous basic invariants, with
\(
\deg(\Phi_i)=d_i
\)
for \(1\le i\le r\). Let \(\xi\) be a primitive \(m\)-th root of unity and, for 
\(1\le i\le r\), let \(a_i\in\ZZ/m\ZZ\) satisfy \(\theta\cdot\Phi_i=\xi^{a_i}\Phi_i\). 
The non-zero restrictions \(\Phi_i|_V\) with
\(
1\le i\le r,
\) and
\(d_i +a_i\equiv0\pmod m
\)
generate \(\kk[V]^G\), whereas
\[
\eta_i:=(\nabla\Phi_i)|_V,
\qquad i \in I_2 = \{j \in \{1,\ldots,r\} \mid 
d_j+a_j\equiv2\pmod m \},
\]
are invariant derivations.
To write the Saito matrix, we define the infinitesimal action. Let
\[
\psi:\fg\otimes_\kk\bU(-1)\to\Der_\bU
\]
be defined by \(\bU\)-linearity and by \(\psi(\bx)(f)(v)
=
d_vf(\bx\cdot v)\),
for all \(\bx\in\fg\), \(f\in\bU\), and \(v\in V\).
Set \(\Derrho=\im(\psi)\), and let \(M_\psi\) denote the coefficient
matrix of \(\psi\); see \S\ref{section:variant}.
Then our application to cyclic gradings is the following result, see Theorem~\ref{thm-periodic-saito}.

\begin{theorem}
Assume that the grading is stable and that \(G\) is semisimple. Then
\[
\Der_\bU(-\log\Delta)
=
\Derrho\oplus
\bigoplus_{i \in I_2}\bU\eta_i.
\]
Let \(M_\Phi\) be the coefficient matrix of the derivations
\((\eta_i\mid i \in I_2)\). Then
\[
M_\Delta=
\begin{pmatrix}
M_\psi\\
M_\Phi
\end{pmatrix}
\quad\text{satisfies}\quad
\det(M_\Delta)=c\,\Delta, \quad \text{for some \(c\in\kk^\times\)}.
\]
\end{theorem}

The theorem applies to several well-studied invariant polynomials,
including the equations of degrees \(120\) and \(126\) defining the
projective duals of \(G(3,9)\) and \(G(4,8)\), respectively; see
\cite{holweck-oeding}. 
However, our method gives a systematic way to compute discriminants of stable polar representations, together with closed determinantal formulas for them, thus extending the reach of explicit computation well beyond these examples.
A particularly striking case is the degree-\(240\) discriminant of
the half-spin representation \(\Sigma_+\) of \(\Spin(16)\), defining the projective
dual of the spinor variety
\(\operatorname{OG}_{+}(8,16)\subset\PP(\Sigma_+)\).
It arises from the stable \(2\)-grading
\(
\fe_8=\fso(V_{16})\oplus\Sigma_+.
\)
Our construction yields a \(128\times128\) Saito matrix whose
determinant is this degree-\(240\) discriminant.
In \S\ref{sec-stable-periodic} we treat further free divisors arising
from stable gradings, notably gradings of \(\fe_6\), \(\fe_7\) and \(\fe_8\) of
orders \(2\), \(3\), and \(5\).

The structure of the paper is as follows. In
\S\ref{section:basic} we fix the setting and recall the
necessary results on derivations, group actions, polar
representations, and complex reflection groups. In
\S\ref{sec-actionderivation} we study cofree representations
and invariant derivations, after introducing variant derivations and
derivations along the quotient fibres in
\S\ref{section:variant}. In
\S\ref{sec-polardiscriminant} we prove the direct-sum
decomposition, Theorem~\ref{thm_direct_sum_dec}, together with the
needed results on invariant derivations. Finally, in
\S\ref{applications} we treat full adjoint discriminants,
actions with finite generic stabilizer, stable periodic gradings, and
the exceptional examples.

We are indebted to M. Yoshinaga for sharing his ideas on the material discussed in this paper.

\section{Preliminary results}

\label{section:basic}

In this section, \(G\) is a connected reductive algebraic group over \(\kk\).

\subsection{Derivations}
\label{sec_der_tensor}

Let \(\bA\) be a finitely generated commutative \(\kk\)-algebra and
let \(\bM\) be an \(\bA\)-module. We write the \(\bA\)-module of \(\bM\)-valued derivations as
\[
\Der_{\bA,\bM}
=
\left\{
\vartheta\in\Hom_\kk(\bA,\bM)
\mid
\vartheta(ab)=a\vartheta(b)+b\vartheta(a)
\text{ for all }a,b\in\bA
\right\}.
\]
There is a natural
identification
\[
\Der_{\bA,\bM}
\simeq
\Hom_\bA(\Omega_{\bA/\kk},\bM).
\]
Since \(\bA\) is finitely generated over \(\kk\), the module
\(\Omega_{\bA/\kk}\) of Kähler differentials is finitely presented. Hence, if \(\bM\) is flat
over \(\bA\), then
\[
\Der_{\bA,\bM}
\simeq
\Der_\bA\otimes_\bA\bM.
\]
Let \(V\) be a finite-dimensional \(\kk\)-vector space and set
\(\bU=\kk[V]\). For \(F\in\bU\), \(w\in V\) and \(v\in V\), we denote by
\(
d_wF\in V^\vee
\)
the differential of \(F\) at \(w\), and by
\(
\partial_vF\in\bU
\)
the directional derivative of \(F\) in the direction \(v\). Thus
\(
\partial_vF(w)=d_wF(v),
\)
for all \(v,w \in V\).
Equivalently, every polynomial \(F\) determines the \textit{directional derivative} map 
\[
D_F:V\to\bU,
\qquad
v\mapsto\partial_vF.
\]
If \(F\) is homogeneous of degree \(d\), then 
\(
D_F:V\to\bU_{d-1}.
\)

\subsubsection{Logarithmic derivations}

Assume now that \(\dim_{\kk}(V)=\ell\), and choose a basis
\((v_1,\ldots,v_\ell)\) with dual coordinates
\((x_1,\ldots,x_\ell)\).
Given a non-zero homogeneous polynomial \(F\in\bU\) of degree \(d>0\), we define the
\(\bU\)-module of \emph{logarithmic derivations of \(F\)}, or derivations
\emph{tangent to the hypersurface defined by \(F\)}:
\begin{align*}
  \Der_\bU(-\log F) &= \left\{\eta \in \Der_\bU \: | \: \eta(F) \in (F)\coloneqq \bU F\right\}.
\end{align*}
We also use
\(\Der_\bU(-\log F)_0 = \left\{\eta \in \Der_\bU \: | \: \eta(F) = 0 \right\}\).
Since we are in characteristic zero, the Euler derivation \(\epsilon = \sum_{i =1}^\ell x_i\frac{\partial}{\partial x_i}\) provides a splitting:
\[
\Der_\bU(-\log F) = \Der_\bU(-\log F)_0 \oplus \bU \epsilon.
\]
Analogously, given a subring \(\bR\) of \(\bU\) and \(f\in \bR\) we define
\[
\Der_{\bR,\bU}(-\log f) = \left\{\eta \in \Der_{\bR,\bU} \: | \: \eta(f) \in (f) \subset \bU \right\}.
\]

\subsubsection{Non-reduced polynomials}

Let \(F_1,\ldots,F_s\in\bU\) be pairwise non-associated irreducible
polynomials. For every \(s\)-tuple
\(a=(a_1,\ldots,a_s)\) of positive integers, set
\[
F^a=F_1^{a_1}\cdots F_s^{a_s},
\qquad
F=F_1\cdots F_s.
\]
Thus \(F\) is the reduced equation associated with \(F^a\). The next lemma
shows that their modules of logarithmic derivations agree. Fix
\(\vartheta\in\Der_\bU\).
\begin{lemma} \label{powers}
  The following conditions are equivalent.
  \begin{enumerate}[label=\roman*)]
  \item \label{powers1} for all \(i \in \{1,\ldots,s\}\), we have \(\vartheta(F_i) \in (F_i)\);
  \item \label{powers2} for all \(s\)-tuples \(a\), we have \(\vartheta(F^a) \in (F^a)\);
  \item \label{powers3} there is an \(s\)-tuple \(a\) such that \(\vartheta(F^a) \in (F^a)\).
  \end{enumerate}
\end{lemma}
\begin{proof}
By the Leibniz rule, \ref{powers1} \(\Rightarrow\) \ref{powers2}, while \ref{powers2} obviously implies \ref{powers3}.
The statement is clear if \(s=1\) and, by iteration, it suffices to show \ref{powers3} \(\Rightarrow\) \ref{powers1} when \(s=2\). Hence, let \(F,G \in \bU\) be coprime polynomials and choose integers \(a,b > 0\); we assume \(\vartheta(F^aG^b) \in (F^a G^b)\) and show \(\vartheta(F) \in (F)\). By assumption, \(F^a\) divides  \(\vartheta(F^aG^b)=aF^{a-1} G^b\vartheta(F) + b F^aG^{b-1}\vartheta(G)\), so \(F^a\) divides \( aF^{a-1}G^b\vartheta(F)\). 
Since \(\kk\) has characteristic zero and \(F\) is coprime to \(G\), we get that \(F\) divides \(\vartheta(F)\).
By symmetry, \(\vartheta(G) \in (G)\) and the lemma is proved.
\end{proof}

\subsubsection{Group action on derivations}\label{sec-actiondev}

Let \(V\) be a finite-dimensional \(G\)-representation and set
\(\bU=\kk[V]\). There is a natural identification
\(
\Der_\bU\simeq \bU\otimes_\kk V,
\)
which is \(G\)-equivariant for the action on derivations, given by
\[
(g\cdot\eta)(f)=g\cdot\eta(g^{-1}\cdot f),
\]
for all \(g\in G\), \(\eta \in \Der_\bU\) and \(f \in \bU\).
In particular, we have
\begin{equation} \label{DerG}
\Der_\bU^G\simeq(\bU\otimes_\kk V)^G.
\end{equation}

If \((v_1,\ldots,v_\ell)\) is a basis of \(V\) with dual coordinates
\((x_1,\ldots,x_\ell)\), then \(v_i\) corresponds to
\(\partial/\partial x_i\).

\begin{lemma}
\label{lem_action_derivatives_inv}
Let \(0\ne F\in\bU_d^G\), with \(d>0\). The directional derivative map
\(
D_F:V\to\bU_{d-1},
\)
is \(G\)-equivariant. Its image is the span of the partial derivatives
of \(F\) in any system of linear coordinates.
If \(V\) is irreducible, then \(D_F\) is injective, and hence
\(\im(D_F)\simeq V\) as a \(G\)-module.
\end{lemma}

\begin{proof}
For \(g\in G\) and \(v,x\in V\), using the \(G\)-invariance of \(F\), we have
\[
\begin{aligned}
(g\cdot\partial_vF)(x)
&=
d_{g^{-1}x}F(v)
=
d_x F(gv)
=
\partial_{gv}F(x).
\end{aligned}
\]
Thus \(D_F\) is \(G\)-equivariant. Its image is the span of the partial
derivatives since, for a basis \((v_1,\ldots,v_\ell)\) with dual
coordinates \((x_1,\ldots,x_\ell)\),
one has
\[
D_F(v_i)=\frac{\partial F}{\partial x_i},
\qquad
1\le i\le\ell.
\]

If \(V\) is irreducible, then \(\ker(D_F)\) is either \(0\) or \(V\).
The second possibility would imply that all directional derivatives
of \(F\) vanish, hence that \(F\) is constant, contrary to \(d>0\).
\end{proof}

\subsection{Polar representations}

\label{section:polar}

We recall the facts about polar representations that will be used
throughout the paper.

\subsubsection{Cartan subspaces}
\label{section:cartan}

Let \(V\) be a finite-dimensional \(G\)-representation, and set
\(\ell=\dim_{\kk}(V)\). An element \(v\in V\) is \emph{semisimple} if its
\(G\)-orbit is closed. 
We let \(\bU=\kk[V]\) and \(V\sslash G=\Spec(\bU^G)\).

\begin{definition}
For \(v \in V\), set
\(
\fh_v
=
\{x\in V\mid\fg\cdot x\subseteq\fg\cdot v\}.
\)
The representation \(V\) is \emph{polar} if there is a semisimple
element \(v\in V\) such that
\begin{equation}
\label{eq_polar_cond}
\dim(V\sslash G)=\dim_{\kk}(\fh_v).
\end{equation}
A subspace \(\fh_v\) satisfying \eqref{eq_polar_cond} is called a
\emph{Cartan subspace}. All Cartan subspaces are \(G\)-conjugate, and
we usually denote one of them simply by \(\fh\). We set
\(
n=\dim_{\kk}(\fh).
\)
\end{definition}
By \cite[Corollary 2.5]{dadok-kac}, if \(V\) is polar, then there is a decomposition
\[
V=\fh\oplus\fg\cdot\fh\oplus\fu,
\]
where \(\fu\) is a \(G_v\)-representation such that
\(\dim(\fu\sslash G_v)=0\).

\begin{remark}
\label{rmk-stability}
We use \emph{stable} in the sense of Dadok--Kac: a representation
of a reductive group is stable if a general orbit is closed,
equivalently, if there exists a closed orbit of maximal dimension.
For a polar representation, in the decomposition
\(
V=\fh\oplus\fg\cdot\fh\oplus\fu
\)
of \cite[Corollary 2.5]{dadok-kac}, stability is equivalent to
\(\fu=0\). Thus a stable polar representation satisfies
\(
V=\fh\oplus\fg\cdot\fh.
\)
A representation is called \emph{locally free} if the stabilizer of a
general point is finite. We reserve the term \emph{GIT-stable} for the
combination of these two conditions, i.e. for the existence of a
closed orbit with finite stabilizer.
An element \(x\in V\) is said to be \textit{stable} if the \(G\)-orbit is closed and its stabiliser in \(G\) is finite.
\end{remark}

We note that, in the terminology used by \cite{panyushev:regular-elements}, representations \textit{with a Cartan subspace} are in fact stable polar representations, see \cite[Corollary 2.5]{dadok-kac}.

\subsubsection{Weyl groups of polar representations}
Let \(V\) be a polar representation of a reductive algebraic group \(G\) and choose a Cartan subspace \(\fh\subset V\). 
We consider the normalizer and centralizer (equivalently, the
pointwise stabilizer)
\[
N_G(\fh)
=
\{g\in G\mid g\cdot\fh=\fh\},
\qquad
Z_G(\fh)
=
\{g\in G\mid g\cdot v=v\text{ for all }v\in\fh\}.
\]
Then the Weyl group is defined as 
\[
W \coloneqq N_G(\fh)/Z_G(\fh).
\]
The Weyl group acts linearly on \(\fh\). 
We will use throughout the following notation:
\[
\bS \coloneqq \kk[\fh].
\]
By \cite[Theorem 2.9]{dadok-kac}, we have the fundamental isomorphism
\[
\bR \coloneqq \bU^G\simeq \bS^W.
\] 
Thus \(V\sslash G\simeq\fh/W\). 

\begin{remark}
\label{rem_stabilizer}
Fix a Cartan subspace \(\fh=\fh_v\) for some semisimple element
\(v\in V\), and let \(w\in\fh\) be general. By
\cite[Lemmas 2.1 and 2.11]{dadok-kac}, the element \(w\) is
semisimple and
\(
\fg\cdot w=\fg\cdot\fh.
\)
On the other hand, the definition of \(\fh=\fh_v\) gives
\(
\fg\cdot\fh=\fg\cdot v.
\)
Therefore
\(
\fh_w=\fh_v=\fh.
\)

Now, if \(g\in G_w\), then
\(
g\cdot\fh_w=\fh_{g\cdot w}=\fh_w,
\)
and consequently
\(
G_w\subseteq N_G(\fh).
\)
The image of \(G_w\) in
\(
W=N_G(\fh)/Z_G(\fh)
\)
fixes \(w\). Since \(W\) is finite and acts faithfully on \(\fh\), the
stabilizer of a general element of \(\fh\) is trivial. Thus
\(
G_w\subseteq Z_G(\fh).
\)
The reverse inclusion is immediate since \(w\in\fh\), and hence
\[
G_w=Z_G(\fh).
\]
\end{remark}

\subsubsection{Regular elements}
\label{sec_reg_elements}
An element \(v\in V\) is \emph{regular} if its orbit has maximal
dimension among the \(G\)-orbits in \(V\).
This convention allows regular elements that are not semisimple; for
example, every non-zero nilpotent element of \(\fsl_2\) is regular.
In \cite{dadok-kac}, the term regular is reserved for
maximal-dimensional semisimple orbits.

The remaining elements are called \emph{irregular}. The next lemma
gives a determinantal description of the irregular locus. Set
\[
s_V\coloneqq \max_{v\in V}\dim(G\cdot v),
\]
so \(s_V\) is the generic orbit dimension of the representation \(V\).
Equivalently, for a general \(v\in V\),
\(
s_V=\dim_{\kk}(\fg)-\dim_{\kk}(\fg_v).
\)
For \(v\in V\), we consider the linear map
\(
\psi_v:\fg\to V,
\)
defined by
\(\psi_v(\bx)=\bx\cdot v,
\)
for all \(\bx \in \fg\).
This agrees with evaluation at \(v\) of the map \(\psi\) defined in the introduction, see also \S\ref{subsubsection:variant}.

\begin{lemma}
\label{lem-irregular-determinantal}
An element \(v\) is regular if and only if
\(
\rk(\psi_v)=s_V.
\)
If \(V\) is a stable polar representation and \(n=\dim_{\kk}(\fh)\), then
\(
s_V=\ell-n.
\)
\end{lemma}

\begin{proof}
For every \(v\in V\), the image of \(\psi_v\) is the tangent space
\(\fg\cdot v\) to the orbit \(G\cdot v\) at \(v\). Hence
\[
\rk(\psi_v)=\dim(G\cdot v),
\]
which proves the first assertion.
Now assume that \(V\) is stable polar. For a general element
\(v\in\fh\), one has
\(
\fg\cdot v=\fg\cdot\fh.
\)
Since stability gives
\(
V=\fh\oplus\fg\cdot\fh,
\)
we obtain
\[
\dim(G\cdot v)=\dim_{\kk}(\fg\cdot\fh)=\ell-n.
\]
Since a general element \(v\) is regular, this is the maximal orbit
dimension, and hence
\(
s_V=\ell-n.
\)
\end{proof}

After fixing bases of \(\fg\) and \(V\), the maps \(\psi_v\) are
represented by a matrix of linear forms on \(V\). By the lemma, the
irregular locus is the determinantal closed subset defined by its
\(s_V\times s_V\) minors.
The next lemma describes the image of the irregular locus in
\(V\sslash G\).

\begin{lemma}
\label{lem_image_hyperplanes}
Let \(V\) be a stable polar representation. There exists a dense open
subset \(\Omega\subset V\sslash G\) such that, for every \(y\in\Omega\),
the fibre \(\chiV^{-1}(y)\) is a single closed regular orbit.
Moreover, the image of the irregular locus is contained in a hypersurface of \(V\sslash G\).
\end{lemma}

\begin{proof}
Consider the finite quotient morphism
\(
q:\fh\to\fh/W\simeq V\sslash G.
\)
By \cite[Lemma~2.11]{dadok-kac}, there is a dense \(W\)-invariant open
subset \(\fh^{\rm reg}\subset\fh\) consisting of regular semisimple
elements. Since \(q\) is finite,
\(
\Omega:=(\fh/W)\setminus q(\fh\setminus\fh^{\rm reg})
\)
is a dense open subset of \(\fh/W\).

Let \(y\in\Omega\) and choose \(v\in q^{-1}(y)\). Then \(v\) is regular
semisimple, so \(G\cdot v\) is a closed orbit of maximal dimension.
Recall that every fibre of the quotient morphism contains a unique
closed orbit, and that the closure of every orbit in the fibre
contains this closed orbit; see
\cite[Theorem~4.2]{brion-etall:luna}.
Now let \(x\in\chiV^{-1}(y)\). If
\(
G\cdot x\neq G\cdot v,
\)
then \(G\cdot v\) lies in the boundary of \(G\cdot x\), and therefore
\(
\dim(G\cdot x)>\dim(G\cdot v),
\)
contradicting the regularity of \(v\). Hence
\[
\chiV^{-1}(y)=G\cdot v.
\]

It follows that the irregular locus does not meet
\(\chiV^{-1}(\Omega)\), so its image is contained in the proper closed
subset
\(
(V\sslash G)\setminus\Omega.
\)
Since \(V\sslash G\simeq\fh/W\simeq\AA^n_{\kk}\), this proper closed subset is
contained in a hypersurface.
\end{proof}

\subsubsection{Visible representations}

We use the convention that visibility is a condition on all fibres of
the quotient morphism.
\begin{definition}
A representation \(V\) of a reductive group \(G\) is called \emph{visible}
if every fibre of
\[
\chiV:V\to V\sslash G
\]
contains finitely many \(G\)-orbits.
\end{definition}

For linear representations of reductive groups, it is enough to check
this condition on the null-cone \(\chiV^{-1}(0)\); see the discussion in
\cite[\S4]{bulois-lehn-lehn-terpereau:symplectic-chevalley}. Thus our
definition agrees with the older convention in which visibility is
formulated only in terms of the null-cone.

The adjoint representations and, more generally, the
theta-representations considered below are visible.

\begin{lemma}
\label{codimension}
Let \(V\) be a visible stable polar representation.
Then the irregular locus has codimension at least two in \(V\).
\end{lemma}

\begin{proof}
Let
\(
\Omega\subset V\sslash G
\)
be the dense open subset provided by
Lemma~\ref{lem_image_hyperplanes}, and set
\(
B=(V\sslash G)\setminus\Omega.
\)
The irregular locus is contained in \(\chiV^{-1}(B)\), and
\(
\dim(B)\leq n-1.
\)

By Lemma~\ref{lem-irregular-determinantal}, the maximal dimension of
a \(G\)-orbit is
\(
\ell-n.
\)
Hence every irregular orbit has dimension at most
\(
\ell-n-1.
\)
Since \(V\) is visible, every fibre of \(\chiV\) contains
only finitely many \(G\)-orbits. Therefore, for every \(y\in B\), the
irregular part of the fibre \(\chiV^{-1}(y)\) has dimension at most
\(\ell-n-1\).
Thus, by the fibre-dimension inequality, the irregular
locus has dimension at most
\[
\dim(B)+(\ell-n-1)
\leq
(n-1)+(\ell-n-1)
=
\ell-2.
\]
Thus its codimension in \(V\) is at least two.
\end{proof}

\subsection{Complex reflection groups}
\label{sec:exponents}

Our reference for complex reflection groups is
\cite{lehrer-taylor:unitary}. The Weyl group
\(W\subset\GL(\fh)\) is generated by pseudo-reflections, namely
finite-order linear transformations whose fixed space is a
hyperplane. By the Chevalley--Shephard--Todd theorem, its invariant
ring is a polynomial ring
\[
\bR=\bS^W\simeq\kk[f_1,\ldots,f_n],
\]
where the basic invariants \(f_1,\ldots,f_n\) may be chosen homogeneous,
of degrees \(d_1\leq\cdots\leq d_n\).

Let \(\bR_+\bS\) be the ideal generated by the positive-degree
\(W\)-invariants. The coinvariant algebra
\[
\bH=\bS/\bR_+\bS
\]
inherits the grading of \(\bS\) and is isomorphic, as a
\(W\)-representation, to the regular representation \(\kk[W]\); see
\cite[\S3]{lehrer-taylor:unitary}. Hence, for every finite-dimensional \(W\)-module \(\bM\), the space
\[
\Hom_{\Wmod}(\bM,\bH)
=
\bigoplus_{d\geq 0}\Hom_{\Wmod}(\bM,\bH_d)
\]
is naturally graded.
Since \(\bH\) is isomorphic to the regular representation of \(W\),
we have
\[
\dim_{\kk}\bigl(\Hom_{\Wmod}(\bM,\bH)\bigr)
=
\dim_{\kk}(\bM).
\]
In particular,
\[
\dim_{\kk}\bigl(\Hom_{\Wmod}(\fh,\bH)\bigr)
=
\dim_{\kk}\bigl(\Hom_{\Wmod}(\fh^\vee,\bH)\bigr)
= 
n.
\]
Choose homogeneous bases
\((u_i\mid 1\le i\le n)\)
and \((w_i\mid 1\le i\le n)\)
of these two spaces, respectively. Up to reordering, for
\(1\le i\le n\), write
\[
e_i=\deg(u_i),\qquad
c_i=\deg(w_i),
\]
with
\(
e_1\leq\cdots\leq e_n\) and
\(c_1\leq\cdots\leq c_n.
\)
The integers \(e_1,\ldots,e_n\) and \(c_1,\ldots,c_n\) are called respectively the
\emph{exponents} and \emph{coexponents} of the reflection representation.
For a polar representation \(V\) with Weyl group \(W\), we write
\[
\Exp(V)=(e_1,\ldots,e_n),
\qquad
\Coexp(V)=(c_1,\ldots,c_n).
\]
In fact, the exponents \(e_1,\ldots,e_n\) are given by the formula
\(e_i=d_i-1\) for \(1\le i\le n\). When the reflection representation is the complexification of a real
one, the \(W\)-modules \(\fh\) and \(\fh^\vee\) are isomorphic, so
exponents and coexponents coincide. Moreover, if \(\fh^W=0\), then
\(c_1=1\): there is no coexponent of degree zero, while
\(\bH_1=\fh^\vee\).

\subsection{Restriction to a Cartan subspace}

\label{section:derivations}

Let \(V\) be a polar representation and choose a Cartan subspace
\(\fh\subset V\). Choose a basis \((v_1,\ldots,v_\ell)\) adapted to
the decomposition
\[
V=\fh\oplus\fg\cdot\fh\oplus\fu,
\]
with \((v_1,\ldots,v_n)\) a basis of \(\fh\), and let
\((x_1,\ldots,x_\ell)\) be the dual coordinates. We call these
\emph{adapted coordinates}. In particular,
\[
x_i|_\fh=0
\qquad\text{for }i>n,
\]
and the restrictions of \(x_1,\ldots,x_n\) form coordinates on \(\fh\).
The following observation will be used when restricting invariant
derivations to \(\fh\).

\begin{lemma}
\label{lem_der_parallel_orbits}
Let \(F\in\bU^G\). Then, for every \(u\in\fg\cdot\fh\), we have
\(
\left.\partial_uF\right|_\fh=0.
\)
\end{lemma}

\begin{proof}
Since \(F\) is \(G\)-invariant, it is constant along every \(G\)-orbit.
Hence, for every \(w\in V\),
\[
d_wF(\fg\cdot w)=0.
\]
For a general \(w\in\fh\), one has
\(
\fg\cdot w=\fg\cdot\fh.
\)
Therefore, for every \(u\in\fg\cdot\fh\),
\(
d_wF(u)=\partial_uF(w)=0
\)
for general \(w\in\fh\). Since \(\partial_uF|_\fh\) is a polynomial
function on \(\fh\), it vanishes identically.
\end{proof}

For adapted coordinates and \(1 \le i \le \ell\), if the basis
vector \(v_i\) belongs to \(\fg\cdot\fh\), the lemma gives
\[
\left.\frac{\partial F}{\partial x_i}\right|_\fh=0.
\]

\subsection{Vinberg representations}\label{section:theta-representations}
We conclude with the facts about Vinberg theta-representations needed
in the applications. Let \(\LL\) be a connected reductive algebraic
group with Lie algebra \(\fL\), and let
\(\theta\in\operatorname{Aut}(\LL)\) have finite order \(m>1\). Fix a primitive \(m\)-th root of unity \(\xi\), and denote again by \(\theta\) the induced automorphism of \(\fL\). We set
\[
\fL_i:=\{x\in\fL\mid \theta(x)=\xi^i x\},
\qquad i\in\ZZ/m\ZZ.
\]
Then we have an \(m\)-cyclic grading of \(\fL\), as follows:
\[
\fL=\bigoplus_{i\in\ZZ/m\ZZ}\fL_i,
\qquad
[\fL_i,\fL_j]\subseteq\fL_{i+j}, \qquad
\forall i,j\in\ZZ/m\ZZ.
\]
Let
\(
G=(\LL^\theta)^\circ,
\)
where \(\LL^\theta\) is the fixed-point subgroup of \(\theta\) and \((\LL^\theta)^\circ\) is its identity component. Then
\(
\fg:=\operatorname{Lie}(G)=\fL_0,
\)
and the adjoint action of \(G\) preserves every graded piece \(\fL_i\),
for \(i\in\ZZ/m\ZZ\).
The \(G\)-representation \(V:=\fL_1\) is called the \textit{Vinberg theta-representation} associated with \((\LL,\theta)\). We say that the grading is \emph{inner} if \(\theta\) is inner, that is, if \(\theta=\operatorname{Int}(g)\) for some \(g\in\LL\).

Choose a \(\theta\)-invariant non-degenerate \(\LL\)-invariant bilinear form
\[
B:\fL\times\fL\to\kk.
\]
It pairs opposite graded pieces perfectly, and therefore induces \(G\)-equivariant isomorphisms
\[
(\fL_{-i})^\vee\simeq\fL_i
\]
for all \(i\in\ZZ/m\ZZ\). We refer to \cite[\S1.3]{Vinberg}.

\begin{proposition} \label{vinberg-is-polar}
Vinberg theta-representations of a reductive algebraic group are visible polar representations.
\end{proposition}

\begin{proof}
Write \(\zeta(-)\) for the centralizer in \(\fL\). For \(x\in\fL_1\),
the tangent space \(\fg\cdot x\) is the orthogonal complement of
\(\zeta(x)\cap\fL_{-1}\) by \cite[Proposition~5]{Vinberg}. If \(\fh\subset\fL_1\) is a Cartan subspace and \(x\in\fh\) is general, then \(\zeta(x)=\zeta(\fh)\) by \cite[\S3.2]{Vinberg}. Proposition~6 of \textit{loc. cit.} gives
\[
\zeta(\fh)\cap\fL_1=\fh\oplus U,
\qquad
U=[\zeta(\fh),\zeta(\fh)]\cap\fL_1,
\]
with \(U\sslash G_x\) zero-dimensional, which gives polarity by \cite[Corollary~2.5]{dadok-kac}. Visibility follows from finiteness of nilpotent orbits, \cite[Proposition~2]{Vinberg}, together with the equivalence of visibility conventions recalled above.
\end{proof}

\begin{definition} \label{def:stable-grading}
Following \cite[Definition~2.5]{laga-romano}, see also \cite{reeder-levy-yu-gross}, we call the periodic grading \emph{stable} if the Vinberg representation \(G\curvearrowright V\) is GIT-stable. Thus a stable grading is, in the terminology of Remark~\ref{rmk-stability}, both Dadok--Kac stable and locally free.
\end{definition}

We warn the reader that an \(m\)-cyclic grading as above may give rise to a stable polar representation \(V=\fL_1\) of \(G=(L^\theta)^\circ\), without being a stable grading in the sense of Definition \ref{def:stable-grading}. 
For instance, if \(n \ge 2\) is an integer, the \(2\)-cyclic outer grading of \(\fL=\fsl_{2n}\) corresponding to the diagram-folding involution gives \(\fg=\fsp_{2n}\) acting on 
\(V\simeq\wedge^{\langle 2\rangle }V_{2n},
\)
where \(V_{2n}\) is a vector space of dimension \(2n\) equipped with a symplectic \(2\)-form \(\omega\). Let \(\Omega\in\wedge^2V_{2n}\) be the invariant bivector corresponding to \(\omega\) under the symplectic identification \(V_{2n}\simeq V_{2n}^\vee\). Then
\(\wedge^{\langle 2\rangle }V_{2n}\) is the quotient of \(\wedge^{2}V_{2n}\) by \(\kk \, \Omega\). 
This action is stable polar, but has generically infinite stabilisers, so it is not GIT-stable.

\medskip

For the remainder of this subsection, assume that \(\LL\) is
semisimple. Associated with a Cartan subspace \(\fh\subset\fL_1\) is a \(\theta\)-stable torus \(T\subset\LL\) whose Lie algebra \(\ft\) contains \(\fh\) and whose dimension is
\[
\dim(T)=\dim_{\kk}(\fh)\varphi(m),
\]
where \(\varphi\) denotes Euler's function. If \(r\) denotes the rank of \(\LL\), define
\[
\theta\text{-}\mathrm{corank}:=r-\dim_{\kk}(\fh)\varphi(m).
\]

\begin{proposition}
\label{prop_thetacorank_stable}
If the \(\theta\)-corank is zero, then the periodic grading is stable.
\end{proposition}

\begin{proof}
If the \(\theta\)-corank is zero, then \(T\) is a maximal torus. For
\(i\in\ZZ/m\ZZ\), let \(\ft_i=\ft\cap\fL_i\). By
\cite[\S3.1]{Vinberg}, we have
\[
\dim_{\kk}(\ft_i)=
\begin{cases}
\dim_{\kk}(\fh),&\gcd(i,m)=1,\\
0,&\gcd(i,m)>1,
\end{cases}
\]
and \(\zeta(\fh)=\zeta(\ft)\).
In particular \(\ft_1=\fh\) and \(\ft_0=0\). Since \(\LL\) is semisimple and \(T\) is maximal, \(\zeta(\ft)=\ft\). Hence we have 
\[
\zeta(\fh)\cap\fL_1=\fh,
\]
so the Vinberg representation is Dadok--Kac stable. For a general \(x\in\fh\), \cite[\S3.2]{Vinberg} gives \(\zeta(x)=\zeta(\fh)\), and therefore
\[
\operatorname{Lie}(G_x)=\fg\cap\zeta(x)=\fL_0\cap\ft=\ft_0=0.
\]
Thus the generic stabilizer is finite, and the grading is stable.
\end{proof}

A full list of simple Lie algebras with zero \(\theta\)-corank, with
\(m\) prime and semisimple degree-zero group, is given in
\cite[Table~1]{benedetti-manivel:discriminants}.

\section{Variant derivations and vertical derivations}
\label{section:variant}

Let \(G\) be an algebraic group acting linearly on a vector space \(V\) of dimension \(\ell\).
We define a module of derivations arising from the \(G\)-action on
\(V\), called \textit{variant derivations}. These vector fields form
the orbit distribution of the action and annihilate every invariant
function.

\subsection{Derivations along the quotient fibres}\label{sec-derorbits}

Consider the GIT quotient map \(\chiV : V \to V \sslash G\) associated with the injection \(\bR = \bU^G \hookrightarrow \bU\). Then the differential of \(\chiV\) can be considered as a morphism of \(\bU\)-modules:
\[
d \chiV : \Der_\bU \to \Der_{\bR,\bU}.
\]
The \textit{derivations along the quotient fibres}, also called
\textit{vertical derivations}, are defined as
\(\Der_\bU^{\chiV}\coloneqq\ker(d\chiV)\). In other words, they are
precisely the derivations that annihilate every invariant function:
\[
\Der_\bU^{\chiV}
=
\{\vartheta\in\Der_\bU\mid\vartheta(f)=0
\text{ for every }f\in\bR\}.
\]

If \(\bR\) is a polynomial algebra, then \(\Der_{\bR,\bU}\) is free. Hence \(\Der_\bU^{\chiV}\) is the kernel of a map between free modules, and therefore it is reflexive.
For every \(f\in\bR\), the differential of \(\chiV\) restricts to a map
\[
d\chiV_{f}:\Der_\bU(-\log f)
\to\Der_{\bR,\bU}(-\log f).
\]
The inclusion \(\Der_\bU^{\chiV} \subset \Der_{\bU}(-\log f)\) implies that \(\ker(d\chiV_f) \simeq \Der_\bU^{\chiV}\). For every \(f\in\bR\), this gives an exact sequence
\begin{equation}\label{eq-derchilogf}
    0 \rightarrow \Der_\bU^{\chiV} \rightarrow \Der_{\bU}(- \log f) \rightarrow \Der_{\bR,\bU}(-\log f).
\end{equation}

\subsection{Variant derivations}

\label{subsubsection:variant}

Let \(\rho:G\to\GL(V)\) be the representation and consider its
differential
\[
d_{\id_G}\rho:\fg\to\gl(V).
\]
For \(\bx\in\fg\), we denote by \(\psi(\bx)\) the linear vector field on
\(V\) induced by \(d_{\id_G}\rho(\bx)\). Equivalently,
\[
\psi(\bx)(f)(v)
=
d_vf(\bx\cdot v),
\qquad
f\in\bU,\quad v\in V.
\]
Extending by \(\bU\)-linearity gives a morphism
\[
\psi:\fg\otimes_\kk\bU(-1)\to\Der_\bU.
\]

Choose a basis \((v_1,\ldots,v_\ell)\) of \(V\), with dual coordinates
\((x_1,\ldots,x_\ell)\), and write
\(
d_{\id_G}\rho(\bx)(v_i)
=
\sum_{j=1}^{\ell}\bx_{i,j}v_j,
\)
for some coefficient matrix \((\bx_{i,j})_{1 \le i,j \le \ell}\).
Then the explicit formula for \(\psi\) reads
\begin{equation}
\label{psi-coordinates}
\psi(\bx)
=
\sum_{1\leq i,j\leq\ell}
\bx_{i,j}x_i\frac{\partial}{\partial x_j}.
\end{equation}

\begin{definition}
The module of \textit{variant derivations} is \(\Derrho \coloneqq \im(\psi) \subset \Der_\bU\).
\end{definition}

Notice that, by definition, \(\Derrho\) is a torsion-free module since it is contained in the free \(\bU\)-module \(\Der_\bU\).
Recall that \(s_V\) denotes the generic orbit dimension of \(V\).
The next lemma computes the rank of the module of variant derivations in terms of \(s_V\). It also says that variant derivations annihilate every invariant function, namely
\[
\Der_\bU^{\chiV} =\bigcap_{f \in \bU^G} \Der_\bU(-\log f)_0.
\]

\begin{lemma}
\label{inclusion}
The module \(\Derrho\) has rank \(s_V\) and
\(
\Derrho\subset\Der_\bU^{\chiV}.
\)
Equivalently, every derivation \(\vartheta\in\Derrho\) annihilates
the invariant subring \(\bR=\bU^G\).
\end{lemma}

\begin{proof}
For \(v\in V\), evaluation at \(v\) gives a linear map
\(
\psi_v:\fg\rightarrow V,
\)
that sends \(\bx \in \fg\) to \(\bx\cdot v\).
The image of this map is the tangent space \(\fg\cdot v\) to the orbit \(G\cdot v\).
Hence, for general \(v\in V\),
\[
\rk(\psi_v)=\dim(G\cdot v)=s_V.
\]
Therefore
\(
\rk(\Derrho)=s_V.
\)
Now let \(f\in\bR=\bU^G\). Since \(f\) is constant along every
\(G\)-orbit, for every \(v\in V\) and \(\bx\in\fg\) one has
\(
d_vf(\bx\cdot v)=0.
\)
By the definition of \(\psi\), we have
\[
\psi(\bx)(f)(v)=d_vf(\bx\cdot v)=0.
\]
Thus \(\psi(\bx)(f)=0\) for every \(\bx\in\fg\).
So the same holds for every \(\vartheta\in\Derrho\).
Hence
\(
\Derrho\subset\Der_\bU^{\chiV}.
\)
\end{proof}

For polar representations, there is a deeper relation between the
module \(\Derrho\) of variant derivations and the module
\(\Der_\bU^{\chiV}\) of derivations along the quotient fibres.

\begin{proposition}
\label{prop_double_dual}
Let \(V\) be a visible stable polar representation of a semisimple connected algebraic group \(G\). Then the module \(\Der_\bU^{\chiV}\) is the reflexive hull of 
\(
\Derrho
\).
\end{proposition}

\begin{proof}
To prove the result, we sheafify all modules under consideration, which is harmless since we are working over an affine space \(V\), so working with modules or sheaves makes no substantial difference. We denote the sheafification of \(\Derrho\) and of \(\Der_\bU^{\chiV}\), respectively, by \(\cT_V^\rho\) and \(\cT_V^{\chiV}\). We note that the sheaf \(\cT_V^{\chiV}\) is reflexive. We will prove that, along an open subset whose complement is of codimension at least two, these sheaves are isomorphic and locally free. This implies that their reflexive hulls are isomorphic, which gives our statement, as \(\Der_\bU^{\chiV}\) is reflexive.

We consider the subset \(R \subset V\) consisting of regular elements of \(V\). By \cite{panyushev:regular-elements}, the map \(d\chiV\) induces a morphism of sheaves, whose restriction to \(R\) is surjective. 
So its kernel \((\cT_V^{\chiV})|_{R}\)
is a vector subbundle of the trivial bundle \(V\otimes \cO_R\), seen as the tangent bundle of \(R\).
Its rank is:
\[
\rk((\cT_V^{\chiV})|_{R})=\ell - \dim(V\sslash G)=\ell-\dim_{\kk}(\fh).
\]
For any \(v \in R\), the fibre \((\cT_V^{\chiV})|_{v}\) is a vector subspace of \(V\), of dimension \(\ell-\dim_{\kk}(\fh)\).
On the other hand, for any \(v\in R\), the fibre of
\(\cT_V^\rho\) at \(v\) is
\(
(\cT_V^\rho)_v=\fg\cdot v.
\)
Since \(v\) is regular, Lemma~\ref{lem-irregular-determinantal} gives
\[
\dim_{\kk}(\fg\cdot v)=s_V=\ell-n.
\]

Thus \((\cT_V^\rho)|_R\) is a vector subbundle of
\(V\otimes\cO_R\) of rank \(\ell-n\), hence of the same rank as
\((\cT_V^{\chiV})|_R\).
After sheafification, Lemma~\ref{inclusion} gives an inclusion
\(\cT_V^\rho
\subset
\cT_V^{\chiV}\).
Along the fibre of any point \(v \in R\), this inclusion gives an injective linear map \((\cT^\rho_V) |_v \to (\cT^{\chiV})|_v\) of vector subspaces of \(V\). Since these subspaces have the same dimension, this inclusion is an isomorphism.
Finally, by Lemma~\ref{codimension}, the complement of \(R\), namely the locus of irregular elements, has codimension at least \(2\). The statement follows.
\end{proof}

\begin{example} 
For reductive groups, the equality in Proposition \ref{prop_double_dual} can fail without stability. As a simple example, we analyse a non-trivial action of \(G=\GG_m\) on \(V=\kk^\ell\), with \(\ell > 0\), via a diagonal representation with weights
\(
(a_1,\ldots,a_\ell)
\).
Write \(
V=V_+\oplus V_-\oplus V_0,
\)
where \(V_+\), \(V_-\), and \(V_0\) are respectively the sums of the positive-, negative-, and zero-weight spaces. Denote their dimensions by \(\ell_+\), \(\ell_-\), and \(\ell_0\), so that
\(
\ell_++\ell_-+\ell_0=\ell.
\)
Then \(V\) is stable in the sense of GIT if 
and only if \(\ell_+>0\) and \(\ell_->0\).
Indeed, decomposing \(v\in V\) as \(v=v_++v_-+v_0\) according to \(
V=V_+\oplus V_-\oplus V_0
\), we have that \(v\) is stable if and only if \(v_+ \ne 0 \ne v_-\).

The Lie algebra \(\fg\) of \(\GG_m\) is \(\kk\), acting diagonally with non-zero weights on 
\(
V_+
\)
and 
\(V_-\).
 If \(\ell_+>0\) and \(\ell_->0\), then \(\dim(V\sslash \GG_m)=\ell-1\).
Choose non-zero vectors \(v_+ \in V_+\) and \(v_-\in V_-\). Let \(v=v_++v_-\) and \(\fh = \kk v + V_0\).
Then \(v\) is semisimple and \(\fh=\{x \in V \mid \fg \cdot x \subset \fg \cdot v\}\), because the action of \(\fg\) on \(V_+\) and \(V_-\) is invertible.
Since \(\dim(\fh)=\ell_0+1\), then \(V\) is polar as soon as \(\ell-\ell_0=2\), which means \(\ell_+=1=\ell_-\).
Conversely, if \(V\) is stable polar, all Cartan subspaces have the same dimension, 
which implies \(\ell_+=1=\ell_-\) by the above argument.

In addition, if \(\ell_+=0\) or \(\ell_-=0\), then \(V\) is polar. Indeed, we take \(\fh=V_0\) and \(v\in V_0\), and see that \(\fh=\{x \in V \mid \fg \cdot x \subset \fg \cdot v\}\). Any vector in \(\fh\) is semisimple, $V\sslash G\cong V_0$ and 
\(\dim(V\sslash G)=\ell_0\), which says that \(V\) is polar. However, it is not stable because, in the decomposition $V = \fh \oplus \fg \cdot\fh \oplus \fu$, we have $\fh = V_0$, $\fg \cdot\fh=0$ and $\fu = V_{-} \neq 0$, if $\ell_{+}=0$, or $\fu = V_{+} \neq 0$, if $\ell_{-}=0$.

Let us look at the derivation modules in the situation \(\ell_+\ell_-=0\) and take for instance \(\ell_-=0\). Write \(\bU=\kk[V]=\kk[x_1,\ldots,x_{\ell_+},y_1,\ldots,y_{\ell_0}]\), with variables adapted to the decomposition \(V=V_+\oplus V_0\).
Now, \(\ker(d \chi)\) consists of derivations annihilating all variables \(y_1,\ldots,y_{\ell_0}\). We have $\psi(1)(f)(v)=d_vf(Av)$ with $A=\mathrm{diag}(a_1,\ldots,a_\ell)$,
so $\psi(1)=\varepsilon$, where \(\varepsilon\) is the Euler-type derivation spelled out in \eqref{eq:torus}, and
$\Derrho=\im(\psi)=\varepsilon\bU$. In conclusion:
\begin{equation}\label{eq:torus}
\Der_\bU^{\chiV}=\bigoplus_{1 \le i \le \ell_+}\bU\frac{\partial}{\partial x_i} \supsetneq
\Derrho=\varepsilon\bU, \qquad \text{with} \qquad \varepsilon = \sum_{1 \le i \le\ell_+} a_i x_i \frac{\partial}{\partial x_i}.
\end{equation}

Summarising: if \(\GG_m\) acts non-trivially and \(V\) is polar, we have \((\Derrho)^{\vee\vee}=\Derrho\neq\Der_\bU^{\chiV}\), unless \(\ell_+=\ell_-=1\).
\end{example}

\section{Cofree representations and invariant derivations}
 \label{sec-actionderivation}

In this section we collect some structural results on cofree
representations and their covariants.
For a cofree representation, we
show that the polynomial ring admits a \(G\)-equivariant decomposition
\(
\bU\simeq \bR\otimes_\kk\bE,
\)
and use it to describe modules of covariants. In the stable polar case,
we identify the \(G\)-module \(\bE\) geometrically, which in turn yields an
explicit description of invariant derivations. These results will be
used in \S\ref{sec-polardiscriminant}.

\subsection{Cofree representations}

\label{section:cofree}

A representation \(V\) of an algebraic group \(G\) is \textit{cofree} if \(\bU=\kk[V]\) is a free module over \(\bR = \bU^G\).
Given a cofree representation \(V\), we will usually write \(\bE\) for a \(\kk\)-vector space such that \(\bU = \bR \otimes_\kk \bE\).

\begin{lemma}
  Let \(V\) be a cofree representation of a reductive algebraic group \(G\).
Then there exists a graded \(G\)-submodule \(\bE\subset \bU\) such that
multiplication induces an isomorphism of graded \(G\)-modules
\(
\bR\otimes_{\kk}\bE \to \bU.
\)
Moreover, each homogeneous component \(\bE_d\) is finite-dimensional over
\(\kk\).
\end{lemma}

\begin{proof}
Set
\(
\bR_+ = \bigoplus_{d>0}\bR_d,
\)
and \(
\bF:=\bU/\bR_+\bU.
\)
Since \(\bR_+\bU\) is a homogeneous \(G\)-stable ideal of \(\bU\), the quotient
\(\bF\) inherits a structure of graded \(G\)-module.
Moreover, for every \(d\in\ZZ_{\geq0}\),
the homogeneous component \(\bF_d\) is finite-dimensional over \(\kk\), and we have an exact sequence of finite-dimensional
\(G\)-modules
\[
0 \to (\bR_+\bU)_d
\to \bU_d
\to \bF_d
\to 0.
\]
Since \(G\) is reductive and \(\kk\) has characteristic zero, \(G\) is
linearly reductive. Hence this sequence splits \(G\)-equivariantly.
Therefore, there exists a \(G\)-submodule
\(
\bE_d\subset \bU_d
\)
which maps isomorphically onto \(\bF_d\). Set
\(
\bE:=\bigoplus_{d\geq 0}\bE_d\subset \bU.
\)
Then \(\bE\) is a graded \(G\)-submodule of \(\bU\), and the quotient map
\(\bU\to\bF\) restricts to an isomorphism
\(
\bE \to \bF.
\)

It remains to show that multiplication induces an isomorphism
\(
\mu : \bR\otimes_\kk\bE\to\bU.
\)
We use the graded Nakayama lemma in the following form: if \(\bM\) is a
positively graded \(\bR\)-module and \(\bM=\bR_+\bM\), then \(\bM=0\). Indeed,
otherwise a nonzero homogeneous element of minimal degree cannot belong
to \(\bR_+\bM\).
Now, let \(C=\Coker(\mu)\). Reducing \(\mu\) modulo \(\bR_+\), we obtain the map
\[
(\bR\otimes_\kk\bE)/\bR_+(\bR\otimes_\kk\bE)
\simeq \bE
\to
\bU/\bR_+\bU,
\]
which is an isomorphism by construction of \(\bE\). Hence
\(
C/\bR_+C=0.
\)
By the graded Nakayama lemma, \(C=0\), so \(\mu\) is surjective.
Next, let \(K=\ker(\mu)\). We have an exact sequence
\[
0\to K
\to \bR\otimes_\kk\bE
\overset{\mu}{\to}\bU
\to 0.
\]
Since \(V\) is cofree, \(\bU\) is a free, hence flat, \(\bR\)-module.
Tensoring this sequence over \(\bR\) with
\(
\kk=\bR/\bR_+
\)
therefore gives an exact sequence
\[
0\to K/\bR_+K
\to \bE
\to \bU/\bR_+\bU
\to 0.
\]
The surjection appearing in the above sequence is an isomorphism by
construction. Thus
\[
K/\bR_+K=0.
\]
The graded Nakayama lemma now gives
\(K=0\).
Therefore \(\mu\) is an isomorphism, and hence
\(
\bU\simeq \bR\otimes_\kk\bE
\)
as graded \(G\)-modules.
\end{proof}

\begin{proposition}\label{prop-cofree-covariants}
Let \(V\) be a cofree representation of a reductive group \(G\), and let
\(\bE\subset \bU\) be a graded \(G\)-submodule such that
\(
\bU\simeq \bR\otimes_\kk \bE
\)
as graded \(G\)-modules. Then, for every finite-dimensional \(G\)-module \(\bM\),
there are natural isomorphisms of \(\bR\)-modules
\[
(\bU\otimes_\kk \bM)^G
\simeq
\bR\otimes_\kk(\bE\otimes_\kk \bM)^G
\simeq
\bR\otimes_\kk\Hom_{\Gmod}(\bM^\vee,\bE).
\]
In particular, \((\bU\otimes_\kk \bM)^G\) is a free \(\bR\)-module.
\end{proposition}

\begin{proof}
Since \(G\) acts trivially on \(\bR=\bU^G\), the decomposition 
\(\bU\simeq\bR\otimes_\kk\bE\) gives
\[
(\bU\otimes_\kk \bM)^G
\simeq
\bR\otimes_\kk(\bE\otimes_\kk \bM)^G.
\]
The second isomorphism is the natural identification
\[
(\bE\otimes_\kk \bM)^G
\simeq
\Hom_{\Gmod}(\bM^\vee,\bE).
\]
\end{proof}

\begin{remark}\label{rem-reflection-cofree}
For a finite pseudo-reflection group \(W\) acting on \(\fh\), the preceding
discussion specializes to the classical decomposition
\[
\bS\simeq\bR\otimes_\kk\bH,
\qquad
\bR=\bS^W,
\]
where \(\bH\) may be identified with the coinvariant algebra
\(
\bS/\bR_+\bS.
\)
As a \(W\)-module, \(\bH\) is isomorphic to the regular representation of
\(W\). We refer to \cite{orlik-terao:arrangements}.
\end{remark}

\subsection{Stable polar representations}

For a cofree representation, let \(\bE\subset\bU\) be a graded
\(G\)-submodule such that
\(
\bU\simeq\bR\otimes_\kk\bE.
\)
For a stable polar representation, this \(G\)-module has a 
geometric interpretation in terms of a general fibre of the quotient
morphism.

\begin{theorem}
\label{thm-polarfreemodule}
Let \(V\) be a stable polar representation, and let
\(\bE\subset\bU\) be a graded \(G\)-submodule such that
\(
\bU\simeq\bR\otimes_\kk\bE.
\)
\begin{enumerate}[label=\roman*)]
\item \label{polarfreemodule-i}
Let \(\fh\) be a Cartan subspace and let \(v\in\fh\) be general. Then
there are \(G\)-module isomorphisms
\[
\bE
\simeq
\kk[\chiV^{-1}(\chiV(v))]
\simeq
\kk[G\cdot v]
\simeq
\kk[G/G_v],
\]
where the fibre is understood scheme-theoretically.
\item \label{polarfreemodule-ii}
If \(G\) is semisimple, the same conclusion holds for every regular
semisimple element \(v\in V\).
\end{enumerate}
\end{theorem}

\begin{proof}
Let
\(
y=\chiV(v)\in V\sslash G.
\)
Let us prove \ref{polarfreemodule-i}. By Lemma~\ref{lem_image_hyperplanes},
for a general element \(v\in\fh\),
the set-theoretic fibre over \(y\) is the single closed regular orbit
\[
\chiV^{-1}(y)=G\cdot v.
\]
Since \(V\) and
\(
V\sslash G\simeq\AA^n_{\kk}
\)
are smooth and \(\chiV\) is dominant, generic smoothness shows that,
as soon as \(v\) is general enough, the scheme-theoretic fibre over \(y\) is
smooth, hence reduced. Therefore
\[
\kk[\chiV^{-1}(y)]\simeq\kk[G\cdot v].
\]

Let us look at \ref{polarfreemodule-ii}. Let \(v\in V\) be regular semisimple. Since \(v\) is
semisimple, the orbit \(G\cdot v\) is closed. Pick an element
\(
x\in\chiV^{-1}(y).
\)
The fibre \(\chiV^{-1}(y)\) contains a unique closed \(G\)-orbit, hence
the unique closed orbit contained in \(\overline{G\cdot x}\) is
\(G\cdot v\). If \(x\notin G\cdot v\), then
\[
G\cdot v
\subset
\overline{G\cdot x}\setminus G\cdot x,
\]
and therefore
\(
\dim(G\cdot v)<\dim(G\cdot x).
\)
This is impossible, since \(v\) is regular and the representation \(V\) is stable, so \(G\cdot v\) has
maximal orbit dimension. Hence
\(
\chiV^{-1}(y)=G\cdot v
\)
set-theoretically.

Every point of \(G\cdot v\) is regular. Moreover,
\(V\sslash G\simeq\AA^n_{\kk}\) is smooth. Thus
\cite[Proposition~4.6]{kraft-schwarz:reduced-nullcone}
implies that \(\chiV\) is smooth at every point of \(G\cdot v\).
Consequently the scheme-theoretic fibre over \(y\) is smooth, hence
reduced, and therefore
\(
\kk[\chiV^{-1}(y)]\simeq\kk[G\cdot v],
\)
scheme-theoretically.

Finally, in either case the \(G\)-equivariant \(\bR\)-module isomorphism
\(
\bU\simeq\bR\otimes_\kk\bE
\)
gives, after tensoring over \(\bR\) with the residue field
\(\kappa(y)=\kk\),
\[
\kk[\chiV^{-1}(y)]
=
\bU\otimes_\bR\kappa(y)
\simeq
(\bR\otimes_\kk\bE)\otimes_\bR\kappa(y)
\simeq
\bE.
\]
Since \(G\cdot v\simeq G/G_v\), the assertions follow.
\end{proof}

\begin{corollary}\label{cor-stable-polar-covariants}
Let \(V\) be a stable polar representation, let \(\fh\) be a Cartan
subspace, and let \(v\in\fh\) be general. For every finite-dimensional
\(G\)-module \(\bM\) there is a natural isomorphism of \(\bR\)-modules
\[
(\bU\otimes_\kk \bM)^G
\simeq
\bR\otimes_\kk \bM^{G_v}.
\]
In particular, \((\bU\otimes_\kk \bM)^G\) is a free \(\bR\)-module of rank
\(
\dim_{\kk}(\bM^{G_v}).
\)
\end{corollary}

\begin{proof}
By Proposition~\ref{prop-cofree-covariants}, we have an isomorphism
\(
(\bU\otimes_\kk \bM)^G
\simeq
\bR\otimes_\kk\Hom_{\Gmod}(\bM^\vee,\bE).
\)
By Theorem~\ref{thm-polarfreemodule}, \(\bE\simeq\kk[G/G_v]\).
Hence Frobenius reciprocity gives
\[
\Hom_{\Gmod}(\bM^\vee,\bE)
\simeq
\Hom_{\Gvmod}(\bM^\vee,\kk)
\simeq
\bM^{G_v};
\]
see \cite[Chapter~1, Example~2.3 and Theorem~2.3(2)]
{timashev:homogeneous}.
\end{proof}

\subsection{Derivations tangent to invariant divisors and invariant derivations}

Here we apply the material discussed in the previous subsections to the study of derivations tangent to invariant divisors and of modules of invariant derivations.

\begin{proposition}\label{prop-tensorequality}
If \(V\) is a cofree representation and \(f \in \bR_d\) with \(d>0\), then
\[
\Der_{\bR,\bU}(-\log f) = \Der_\bR(-\log f) \otimes_\bR \bU.
\]
In particular, this holds if \(V\) is a polar representation.
\begin{proof}
Consider the \(\bR\)-linear map
\(
\Der_\bR\to \bR/(f),
\)
defined by sending a derivation
\(\eta
\)
to
\(\eta(f)\bmod f.
\)
Its kernel is \(\Der_\bR(-\log f)\). By assumption, \(\bU\) is a free,
hence flat, \(\bR\)-module. Tensoring the resulting exact sequence with
\(\bU\) therefore identifies
\(
\Der_\bR(-\log f)\otimes_\bR\bU
\)
with the kernel of
\[
\Der_\bR\otimes_\bR\bU
\to
(\bR/(f))\otimes_\bR\bU.
\]
The natural identification from \S\ref{sec_der_tensor} is
\(
\Der_\bR\otimes_\bR\bU
\simeq
\Der_{\bR,\bU}.
\)
Under this identification, the kernel is precisely
\(\Der_{\bR,\bU}(-\log f)\). Finally, polar representations are
cofree by \cite[Theorem~2.10]{dadok-kac}.
\end{proof}
\end{proposition}

\begin{corollary}\label{cor-cofree-invariant-derivations}
Let \(V\) be a cofree representation of a reductive group \(G\), and let
\(\bE\subset\bU\) be as in Proposition~\ref{prop-cofree-covariants}.
Then \(\Der_\bU^G\) is a free \(\bR\)-module, more precisely:
\[
\Der_\bU^G
\simeq
\bR\otimes_\kk
\Hom_{\Gmod}(V^\vee,\bE).
\]
\end{corollary}

\begin{proof}
We have an equivariant identification
\(
\Der_\bU\simeq\bU\otimes_\kk V.
\)
Then, applying Proposition~\ref{prop-cofree-covariants} with \(\bM=V\),
we get 
\(
\Der_\bU^G
\simeq
(\bU\otimes_\kk V)^G
\simeq
\bR\otimes_\kk
\Hom_{\Gmod}(V^\vee,\bE).
\)
\end{proof}

\begin{corollary}\label{cor-freeinvariant}
Let \(V\) be a stable polar representation and let \(v\in\fh\) be general.
Then
\[
\Der_\bU^G
\simeq
\bR\otimes_\kk V^{G_v}.
\]
In particular, \(\Der_\bU^G\) is a free \(\bR\)-module of rank
\(
\dim_{\kk}(V^{G_v}).
\)
\end{corollary}

\begin{proof}
Apply Corollary~\ref{cor-stable-polar-covariants} with \(\bM=V\).
\end{proof}

\begin{remark} \label{rmk-stable-polar-covariants} 
For the action of the pseudo-reflection group \(W\) on \(\fh\), the same
argument gives
\[
\Der_\bS^W
\simeq
(\bS\otimes_\kk\fh)^W
\simeq
\bR\otimes_\kk
\Hom_{\Wmod}(\fh^\vee,\bH).
\]
Since \(\bH\) is the regular representation of \(W\), this is a free
\(\bR\)-module of rank \(\dim_{\kk}(\fh)\). Its homogeneous generator degrees are
the coexponents of \(W\).
\end{remark}

\section{Polar discriminants}
\label{sec-polardiscriminant}

Let \(V\) be a polar \(G\)-representation, let \(\fh\subset V\) be a Cartan
subspace and let \(W\) be its Weyl group. Denote by
\(H_1,\ldots,H_N\) the reflecting hyperplanes of \(W\), choose
\(\delta_j\in\fh^\vee\) with \(H_j=\ker(\delta_j)\), and let \(m_j\) be the
order of the pointwise stabilizer of \(H_j\), for \(1\le j\le N\). Then
\[
\delta=\prod_{j=1}^N\delta_j^{m_j}
\]
is \(W\)-invariant and hence defines an element of
\(
\bR=\bS^W\simeq\bU^G.
\)
We denote by \(\Delta\in\bU^G\) the corresponding invariant polynomial
on \(V\), set
\[
D=\VV(\Delta)\subset V,
\]
and call \(\Delta\) the \emph{polar discriminant equation} and \(D\) the
\emph{polar discriminant divisor}. The element \(\delta\) defines the
branch divisor of \(\fh\to\fh/W\).
Our main result describes the logarithmic derivations along \(\Delta\)
as the direct sum of the derivations along the fibres of the quotient
map and the pull-back of the logarithmic derivations along this branch
divisor.

\begin{theorem}
\label{thm_direct_sum_dec}
Let \(V\) be a stable polar representation of a reductive group \(G\). Then \(\Der_\bU^G\) is a free \(\bR\)-module of rank
\(\dim_{\kk}(V^{G_v})\), where \(v\in\fh\) is general, and  \(\Der_\bR(-\log\delta) \otimes_\bR \bU\) is a free \(\bU\)-module, isomorphic to \(\bigoplus_{c \in \Coexp(V)} \bU(-c)\).
In addition, if the restriction map
\[
\pi_G:\Der_\bU^G \to \Der_\bR(-\log\delta)
\]
is surjective, then there is 
direct-sum decomposition
\[
\Der_\bU(-\log \Delta) \simeq \Der_{\bU}^{{\chiV}} \oplus  \Der_\bR(-\log\delta) \otimes_\bR \bU.
\]
In particular, if
\(\fh = V^{G_v}\),  then $\pi_G$ is an isomorphism.
\end{theorem}

\begin{remark}
Before starting the proof of the above theorem, we show that it applies to orthogonal representations (see Definition \ref{def_ortho} below). In fact, we will show that $\pi_G$ is an isomorphism in the two main applications of the paper, namely adjoint actions of semisimple groups (which are orthogonal because of the existence of the Killing form), and stable gradings (see \S \ref{prop-periodic-coexponents}).
\end{remark}

\begin{definition}
\label{def_ortho}
A polar $G$-representation $V$ is \emph{orthogonal} if there exists a $G$-invariant non-degenerate symmetric form on $V$.
\end{definition}
\begin{proposition}\label{prop-surjadjoint}
Let $V$ be an orthogonal stable polar $G$-representation. Then the restriction map 
\(
\pi_G:\Der_\bU^G \to \Der_\bR(-\log\delta)
\)
is an isomorphism. 
\end{proposition}

\begin{proof}
Since \(V\) is orthogonal, there is a $G$-invariant symmetric bilinear form \(B\) that identifies $V\cong V^\vee$. Moreover, \(B\) induces a $W$-invariant non-degenerate symmetric bilinear form on $\fh$. Thus $\Der_\bU^G\cong (\bU\otimes V)^G\cong (\bU\otimes V^\vee)^G$, which implies that $\Der_\bU^G$ can be identified with differentials of $G$-invariant functions in $\bR$. Analogously, we have $\Der_\bS^W\cong (\bS\otimes \fh)^W\cong (\bS\otimes \fh^\vee)^W$, i.e., $\Der_\bS^W$ can also be identified with differentials of $W$-invariant functions in $\bR$. Using these identifications and Chevalley's Theorem, we see that the restriction map realizes an isomorphism $\Der_\bU^G\xrightarrow{\sim}\Der_\bS^W$. Composing with $\pi_W$, we have that $\pi_G$ is an isomorphism by Theorem \ref{thm_terao_der0}.
\end{proof}

We begin by describing geometric properties of $\Delta$ and its avatar $\delta$.

\begin{remark}
For every \(j\in\{1,\ldots,N\}\), at a general point of \(H_j\), the
stabilizer is cyclic of order \(m_j\),
and the quotient map is transversally of the form
\(
t\longmapsto t^{m_j}.
\)
Thus \(\delta\), as an element of \(\bR\), is a reduced equation of the
branch divisor, whereas its pull-back to \(\fh\) has divisor
\(
\sum_{j=1}^N m_j H_j.
\)
By Lemma~\ref{powers}, these multiplicities do not affect the module
of logarithmic derivations.
\end{remark}

\begin{lemma}\label{lem-polar-discriminant-reduced}
Let \(V\) be a visible stable polar representation of a semisimple
algebraic group \(G\). Then the polar discriminant divisor
\(D=\VV(\Delta)\) is reduced. Equivalently,
\(D=D_\red\) and \(\Delta\) is square-free.
\end{lemma}

\begin{proof}
Let \(B\subset V\sslash G\) be the reduced branch divisor defined by
\(\delta\), and let \(Z\) be an irreducible component of \(D_\red\).
By Lemma~\ref{codimension}, the general point of \(Z\) is regular,
because the irregular locus has codimension at least two. At such a
point, the differential of \(\chiV\) is surjective by
\cite[Corollary 1]{panyushev:regular-elements}; hence \(\chiV\) is smooth there.
The pull-back of a reduced divisor by a smooth morphism is reduced.
It follows that \(\chiV^{-1}(B)=D\) is reduced at the general point of
\(Z\), so \(Z\) occurs in \(D\) with multiplicity one. This holds for
every irreducible component \(Z\), and therefore \(D=D_\red\), or
equivalently \(\Delta\) is square-free.
\end{proof}

\begin{remark}\label{rmk-order6-e8-reducedness}
The semisimplicity assumption in Lemma~\ref{lem-polar-discriminant-reduced}
is essential.
Consider the stable inner \(6\)-grading of \(\fe_8\) corresponding to
the reflection group \(G_{32}\), see
\cite[Tables~5 and~21]{reeder-levy-yu-gross}.
The Weyl group of the associated polar representation is \(G_{32}\).
Its basic invariant degrees are \(12,18,24,30\).
One has
\(\dim_{\kk}(\fg)=40\) and, since the grading is stable of rank \(4\),
\(\dim_{\kk}(V)=44\).
The ambient \(E_8\)-degrees congruent to \(2\) modulo \(6\) are
\(2,8,14,20\), so Proposition~\ref{prop-periodic-gradients} gives four
invariant derivations of degrees \(1,7,13,19\). Together with the
infinitesimal action they give \(44\) generically independent
derivations tangent to \(\Delta\), whose coefficient determinant has
degree
\[
40+(1+7+13+19)=80.
\]
On the other hand, the pull-back \(\Delta\) of the branch equation for
\(G_{32}\) has degree \(120\). Hence \(\Delta\) cannot be reduced: if it
were, the determinant above would have to be divisible by \(\Delta\).
By Lemma~\ref{powers}, however, its multiplicities do not affect
\(\Der_\bU(-\log\Delta)\), which only depends on the reduced divisor
\(D_\red\). This example also explains why the semisimplicity hypothesis
in the Saito-matrix statement below cannot be ignored.
\end{remark}

\subsection{Invariant derivations and the branch divisor}

Restriction of derivations gives an \(\bR\)-linear morphism
\[
\pi_W:\Der_\bS^W\to\Der_\bR.
\]
The following statement summarizes some results from
\cite{orlik-terao:arrangements}. It identifies the logarithmic
derivations along the branch divisor, \(\Der_\bR(-\log\delta)\), with
the image of the invariant derivations on \(\fh\). The reduced inverse
image of this divisor is the reflection arrangement of \(W\).
This is the basic input from the hyperplane arrangement associated with the pseudo-reflection group used
in the proof of Theorem~\ref{thm_direct_sum_dec}.

\begin{theorem}\label{thm_terao_der0}
The restriction morphism induces an isomorphism
\(
\pi_W:
\Der_\bS^W
\to
\Der_\bR(-\log\delta).
\)
These modules are free of rank \(n=\dim_{\kk}(\fh)\). More precisely,
\[
\Der_\bR(-\log\delta)
\simeq
\bigoplus_{c\in\Coexp(V)}\bR(-c).
\]
\end{theorem}

\begin{proof}
The identification by restriction follows from
\cite[Lemma 6.55 and Theorem 6.57]{orlik-terao:arrangements},
and freeness follows from
\cite[Lemma 6.48 and Corollary 6.58]{orlik-terao:arrangements}.
By Remark~\ref{rmk-stable-polar-covariants}, the homogeneous
generator degrees of \(\Der_\bS^W\) are the coexponents of \(W\),
hence the same holds for \(\Der_\bR(-\log\delta)\).
\end{proof}

\subsection{Lifting invariant derivations and the quotient}

Restriction to the invariant ring defines an \(\bR\)-linear map
\[
\pi_G:\Der_\bU^G\to\Der_\bR.
\]
The next lemma identifies the image of \(\pi_G\) with
\(\Der_\bR(-\log\delta)\). In particular, every logarithmic derivation
along the branch divisor admits a \(G\)-invariant lift to \(V\).
This is one of the main ingredients in the proof of
Theorem~\ref{thm_direct_sum_dec}.

\begin{lemma}\label{lemma-surj}
Let \(V\) be a stable polar representation.
Then the image of
\(
\pi_G
\)
is contained in \(\Der_\bR(-\log\delta)\). Moreover, \(\Der_\bU^G\) is a free \(\bR\)-module of rank \(\dim_{\kk}(V^{G_v})\), where \(v\) is a general element of \(\fh\).
\end{lemma}

\begin{proof}
Pick a general element \(v\in\fh\).
By Corollary~\ref{cor-freeinvariant},
\(
\Der_\bU^G\simeq \bR\otimes_\kk V^{G_v}.
\)
This proves the freeness statement and the rank formula.

Since \(V\) is a stable polar representation, we have
\(
V=\fh\oplus\fg\cdot\fh.
\)
Denote by \(p : V \to \fh\) the projection and \(i : \fh \to V\) the injection.
Given an invariant derivation,
\(
\eta \in \Der_\bU^G,
\)
viewing \(\eta\) as a polynomial vector field \(V \to V\), we put \(\Pi(\eta) = p \circ \eta \circ i : \fh \to \fh.\)

We claim that \(\Pi(\eta) \in \Der_\bS^W\).
Indeed, we first observe that both summands of \(V = \fh \oplus \fg \cdot \fh\) are preserved by \(N_G(\fh)\). This is clear for \(\fh\),
while for \(\fg\cdot\fh\) it follows from the equivariance of the
infinitesimal action. Hence the projection
\(
p:V \to \fh
\)
is \(N_G(\fh)\)-equivariant.
Recall that the identification
\(
\Der_\bU\simeq\bU\otimes_\kk V
\)
is \(G\)-equivariant. 
Under this identification, the map \(\Pi : \Der_\bU\to \Der_\bS\) is \(N_G(\fh)\)-equivariant.
Since \(\eta\) is \(G\)-invariant, its image \(\eta_\fh=\Pi(\eta)\) is
\(N_G(\fh)\)-invariant.
Since \(Z_G(\fh)\) acts trivially on \(\fh\), the action of \(N_G(\fh)\)
on \(\bS=\kk[\fh]\), and hence on \(\Der_\bS\), factors through
\(
W=N_G(\fh)/Z_G(\fh).
\)
Therefore
\(
\eta_\fh = \Pi(\eta)\in\Der_\bS^{N_G(\fh)}=\Der_\bS^W.
\)

Next, we note that \(\pi_G = \pi_W \circ \Pi\). Indeed, for all \(F \in \bR\) and \(x \in \fh\), we have 
\(d_x F(\fg \cdot \fh)=0\), see
Lemma~\ref{lem_der_parallel_orbits}. Since \(\eta(x)-p(\eta(x))\in \fg \cdot \fh\), 
we have thus 
\[
\bigl(\pi_G(\eta)(F)\bigr)(x)
=d_xF(\eta(x))
=d_xF\bigl(p(\eta(x))\bigr)
=d_x(F|_\fh)\bigl(\Pi(\eta)(x)\bigr).
\]
This says precisely that \(\pi_G = \pi_W \circ \Pi\).
By Theorem~\ref{thm_terao_der0}, this already shows that
\(
\pi_G(\eta)\in\Der_\bR(-\log\delta).
\)
\end{proof}

\begin{proof}[Proof of Theorem~\ref{thm_direct_sum_dec}]
By Lemma~\ref{lemma-surj}, the surjection
\(
\pi_G:\Der_\bU^G\twoheadrightarrow\Der_\bR(-\log\delta)
\)
admits an \(\bR\)-linear section, since
\(\Der_\bR(-\log\delta)\) is free.
Let
\[
\sigma :\Der_\bR(-\log\delta)\to \Der_\bU^G
\]
be such a section.
Under the identification
\(
\bR=\bS^W\simeq\bU^G,
\)
the elements \(\delta\) and \(\Delta\) correspond to one another. Hence,
 for all \(\vartheta\in\Der_\bR(-\log\delta)\), we have
\[
\sigma(\vartheta)(\Delta)
=
\vartheta(\delta)
\in
\delta\bR
=
\Delta\bR
\subset
\Delta\bU.
\]
Thus
\(
\sigma(\vartheta)\in\Der_\bU(-\log\Delta).
\)
After extending scalars to \(\bU\), we obtain a \(\bU\)-linear section
\[
\Der_\bR(-\log\delta)\otimes_\bR\bU
\to
\Der_\bU(-\log\Delta).
\]
Now the exact sequence~\eqref{eq-derchilogf}, together with
Proposition~\ref{prop-tensorequality}, gives
\[
0\to
\Der_\bU^{\chiV}
\to
\Der_\bU(-\log\Delta)
\to
\Der_\bR(-\log\delta)\otimes_\bR\bU
\to 0,
\]
and the section above splits this sequence. Hence
\[
\Der_\bU(-\log\Delta)
\simeq
\Der_\bU^{\chiV}
\oplus
\left(
\Der_\bR(-\log\delta)\otimes_\bR\bU
\right).
\]

The description of the second summand follows from
Theorem~\ref{thm_terao_der0}. Finally, if \(V^{G_v}=\fh\), then
Corollary~\ref{cor-freeinvariant} and
Theorem~\ref{thm_terao_der0} show that the source and target of
\[
\pi_G:\Der_\bU^G\twoheadrightarrow\Der_\bR(-\log\delta)
\]
are free \(\bR\)-modules of the same rank \(n\). Therefore \(\pi_G\) is an
isomorphism.
\end{proof}

For a hypersurface \(E=\VV(F)\) in a smooth affine variety \(X\), denote by
\(\cT_X\langle E\rangle\) the sheaf associated with the module
\(\Der(-\log F)\) of logarithmic derivations; see \cite{Sernesi}. Let
\(B=\VV(\delta)\subset V\sslash G\) be the branch divisor. Since
logarithmic derivations are insensitive to multiplicities by
Lemma~\ref{powers}, the differential of \(\chiV\) induces a map
\[
\cT_V\langle D_\red\rangle
\to
\chiV^*\cT_{V\sslash G}\langle B\rangle.
\]
Theorem~\ref{thm_direct_sum_dec} has the following sheaf-theoretic
form.

\begin{corollary}\label{main-corollary}
Let \(V\) be a stable polar representation of a reductive algebraic
group \(G\). Then \(\cT_{V\sslash G}\langle B\rangle\) is free,
with homogeneous basis degrees \(\Coexp(V)\), and, if the restriction of invariant derivations is surjective, there is a generally
noncanonical direct-sum decomposition
\[
\cT_V\langle D_\red\rangle
\simeq
\chiV^*\cT_{V\sslash G}\langle B\rangle
\oplus\ker(d\chiV).
\]
\end{corollary}

\section{Applications to adjoint actions, finite stabilisers and stable gradings}

\label{applications}

As mentioned before, it is natural to expect that the module of logarithmic derivations of the discriminant can be decomposed in terms of 
\(G\)-invariant and \(G\)-variant derivations. We here consider two settings, namely adjoint representations and actions with finite generic stabilizers, where this decomposition takes place and affords a full description of the module of derivations tangent to the associated discriminant divisor.

\subsection{Adjoint discriminants}

\label{section:adjoint}

Let \(G\) be a semisimple algebraic group over \(\kk\), acting on its Lie
algebra \(\fg\) by the adjoint representation. Set \(\bU=\kk[\fg]\), let
\(\fh\) be a Cartan subalgebra of \(\fg\), and let \(W\) be the associated
Weyl group. Let \(\Delta\in\bU^G\) be the invariant polynomial
corresponding to the product of all roots of \(\fg\) under the Chevalley
isomorphism
\(
\kk[\fg]^G\simeq\kk[\fh]^W.
\)
We call its zero locus \(\VV(\Delta)\subset\fg\) the \emph{full adjoint
discriminant}.
We refer to \cite[Chapter~3]{slodowy}.

We now describe \(\Der_\bU(-\log\Delta)\) in terms of invariant and
variant derivations, and give a minimal graded free resolution of
\(\Derrhoad\).

\begin{theorem} \label{theorem:adjoint}
Let \(G\) be a semisimple linear algebraic group. Then
\[
\begin{aligned}
\Der_\bU(-\log \Delta)
&\simeq \Derrhoad \oplus \Der_\bU^G\otimes_\bR \bU 
\simeq \Derrhoad \oplus \bigoplus_{e \in \Exp(\fg)}\bU(-e).
\end{aligned}
\]
Moreover, the \(\bU\)-module \(\Derrhoad\) has a minimal graded free resolution of the form:
\[
0 \to \bigoplus_{e \in \Exp(\fg)}\bU(-e-1)
\to \fg \otimes \bU(-1) \to \Derrhoad \to 0.
\]
In particular, \(\Der_\bU(-\log \Delta)\) has projective dimension one.
\end{theorem}

\begin{remark}
For simple groups of simply laced type, this description was proved in
\cite{benedetti-faenzi-marchesi:adjoint-arxiv} by a case-by-case
argument using projective duality. The proof below is uniform and also
applies to semisimple groups of non-simply-laced type.
When \(\fg\) is not simply laced, the full adjoint discriminant considered
here is reducible. The dual of the adjoint variety, studied in
\textit{loc. cit.}, is the irreducible component corresponding to the
\(W\)-orbit of long roots.
\end{remark}

Let \(n=\dim_{\kk}(\fh)\) and choose homogeneous basic invariants of
\(G\), so that
\[
\bR=\bU^G=\kk[F_1,\ldots,F_n].
\]
For \(1\le i\le n\), write \(\deg(F_i)=e_i+1\),
where \(e_1,\ldots,e_n\) are the exponents of \(\fg\). Choose a
non-degenerate invariant symmetric bilinear form \(B\) on \(\fg\), use
it to identify \(\fg\) with \(\fg^\vee\), and define the polynomial map
\(
\nu_i=\nabla F_i:\fg\to\fg
\)
by
\[
B\bigl(\nu_i(x),y\bigr)=d_x F_i(y)
\qquad
\text{for all }x,y\in\fg.
\]
In the adjoint case, \(V^{G_v}=\fh\) for a general
\(v\in\fh\), and the Weyl group is real. Thus
Theorem~\ref{thm_direct_sum_dec} shows that \(\Der_\bU^G\) is free with
homogeneous generator degrees \(e_1,\ldots,e_n\). Kostant's regularity
criterion \cite[Theorem~9]{konstant:polynomialrings} implies that the
gradients \((\nu_i\mid 1\le i\le n)\) are linearly independent at every
regular element. Comparing degrees therefore shows that they form a
homogeneous \(\bR\)-basis of \(\Der_\bU^G\).
Recall from \S\ref{section:variant} that, in the adjoint case,
the map defining variant derivations is
\begin{equation}
  \label{formula-for-ad}
\psi:\fg\otimes\bU(-1)\to\fg\otimes\bU,
\qquad
\psi(\varphi)(x)=[\varphi(x),x].
\end{equation}

\begin{lemma}
\label{lemma-kerad}
For every \(i\in\{1,\ldots,n\}\), the polynomial map \(\nu_i\), regarded
as an element of \((\fg\otimes\bU)(-1)\), belongs to \(\ker(\psi)\).
\end{lemma}

\begin{proof}
Fix \(x,y\in\fg\). The invariance of \(B\) gives
\[
B\bigl([\nu_i(x),x],y\bigr)
=B\bigl(\nu_i(x),[x,y]\bigr)
=d_x F_i([x,y]).
\]
The last term vanishes because \(F_i\) is \(G\)-invariant and \([x,y]\) is
tangent to the adjoint orbit through \(x\). Since \(B\) is
non-degenerate, \([\nu_i(x),x]=0\) for every \(x\in\fg\). This is exactly
the assertion that \(\nu_i\in\ker(\psi)\).
\end{proof}

\begin{proposition}
\label{prop_ex_seq}
For adjoint representations, we have
\[
\ker(\psi)\simeq
\Der_\bU^G\otimes_\bR\bU(-1).
\]
In particular, there is an equivariant exact sequence
\begin{equation}
\label{res-bA}
0\to
\Der_\bU^G\otimes_\bR\bU(-1)
\to\fg\otimes\bU(-1)
\to\Derrhoad\to0.
\end{equation}
\end{proposition}

\begin{proof}
Put \(K = \ker(\psi)\). Then \(K(1)\) is a submodule of \(\fg \otimes \bU\) and, identifying \(\fg \otimes \bU\) with polynomial mappings \(\fg \to \fg\), in view of \eqref{formula-for-ad}, we get:
\[
K(1) =\left\{
\varphi:\fg\to\fg\text{ polynomial}
\ \middle|\ [\varphi(x),x]=0\text{ for every }x\in\fg
\right\}.
\]
Now, we use a result of Dixmier, asserting that this polynomial centralizer module is isomorphic to 
\(
\bigoplus_{i=1}^n\bU\nu_i.
\)
This formulation is recorded explicitly in
\cite[p. 401]{charbonnel-moreau:index-centralizers}, where it is
attributed to \cite{dixmier:adjoint-vector-fields}. Since the same maps
form an \(\bR\)-basis of \(\Der_\bU^G\), we obtain
\[
\ker(\psi)
\simeq \bigoplus_{i=1}^n\nu_i \cdot \bU(-1)
\simeq \Der_\bU^G\otimes_\bR\bU(-1).
\]
The exact sequence follows from the definition
\(\Derrhoad=\im(\psi)\).
\end{proof}

\begin{proof}[Proof of Theorem~\ref{theorem:adjoint}]
Since \(\deg(\nu_i)=e_i\), for \(1 \le i \le n\), Proposition~\ref{prop_ex_seq} gives the
graded free resolution
\[
0\to
\bigoplus_{i=1}^n\bU(-e_i-1)
\to\fg\otimes\bU(-1)
\to\Derrhoad\to0.
\]
It is minimal because every \(e_i\) is positive, so all entries of the
first map belong to the homogeneous maximal ideal of \(\bU\).

It remains to prove that \(\Derrhoad\) is reflexive. It is torsion-free,
since it is a submodule of \(\Der_\bU\). As \(\bU\) is regular, the
criterion of \cite[Proposition~1.3]
{hartshorne:stable-reflexive-sheaves} reduces reflexivity to Serre's
condition \((S_2)\).
Let \(\mathfrak p\subset\bU\) be a prime of height at least three.
Localizing the resolution above and applying the Auslander--Buchsbaum
formula gives
\[
\depth\bigl((\Derrhoad)_\mathfrak p\bigr)
\geq \operatorname{ht}(\mathfrak p)-1
\geq2.
\]
It remains to consider primes of height two. The nonregular locus
\(\fg\setminus\fg_{\mathrm{reg}}\) has codimension three; see
\cite[Theorem~2.14]{moreau:sheets}. Hence every height-two point lies
in \(\fg_{\mathrm{reg}}\). Over the regular locus the infinitesimal
action has constant rank, by Lemma~\ref{lem-irregular-determinantal},
so its image \(\Derrhoad\) is locally free. Thus
\(
\depth\bigl((\Derrhoad)_\mathfrak p\bigr)=2
\)
for every height-two prime \(\mathfrak p\). This proves \((S_2)\), and
hence \(\Derrhoad\) is reflexive.

The adjoint representation is visible, stable, and polar. Proposition
\ref{prop_double_dual} therefore gives
\(
\Derrhoad=\Der_\bU^{\chiV}.
\)
Note that Theorem~\ref{thm_direct_sum_dec}, applies due to Proposition \ref{prop-surjadjoint}. This yields
\[
\Der_\bU(-\log\Delta)
\simeq
\Derrhoad\oplus
\left(\Der_\bU^G\otimes_\bR\bU\right).
\]
Finally,
\(
\Der_\bU^G\simeq\bigoplus_{i=1}^n\bR(-e_i),
\)
which gives the asserted graded decomposition. Since the second
summand is free and \(\Derrhoad\) has projective dimension one, so does
\(\Der_\bU(-\log\Delta)\).
\end{proof}

\subsection{Actions with finite generic stabilizer}

\label{section:finite-stabilizer}

Let \(V\) be a polar representation of a semisimple algebraic group \(G\). Then, we have a Cartan subspace \(\fh \subset V\) and a Weyl group \(W\), which is a finite pseudo-reflection group. We recall from \S\ref{sec:exponents} that the reflecting hyperplanes of \(W\) form a free arrangement, with generators of degrees \((c_1, \ldots, c_n)\). These agree with the coexponents of \(V\) and one writes \(\Coexp(V)=(c_1, \ldots, c_n)\). If the generic stabilizer of the action of \(G\) on \(V\) is finite, the module of \(G\)-variant derivations is also free. 
Indeed, in this case the morphism \(\psi: \fg\otimes \bU(-1) \to V\otimes \bU\) is injective. This gives the next result.

\begin{theorem}
\label{thm_finite_stab}
Let \(V\) be a visible stable polar representation of a semisimple group \(G\) with finite generic stabilizer and surjective restriction of invariant derivations. Then
the derivation module \(\Der_\bU(-\log \Delta)\) is free, more precisely:
\[
\begin{aligned}
\Der_\bU(-\log \Delta)
&\simeq \fg\otimes \bU(-1) \oplus \bigoplus_{c \in \Coexp(V)}\bU(-c).
\end{aligned}
\]
\end{theorem}

\begin{proof}
The result follows directly from Theorem~\ref{thm_direct_sum_dec} once we prove that \(\Der_\bU^{\chiV} \cong \fg\otimes \bU(-1)\). Notice that, under our hypothesis, the map \(\psi : \fg\otimes \bU(-1)\to V\otimes \bU\cong \Der_\bU\) given by the infinitesimal action of \(\fg\) is injective. Indeed, this map is injective if it is generically injective and, in turn, this is guaranteed by the finiteness of the generic stabilizer. This implies that \(\Derrho \cong \fg\otimes \bU(-1)\). Therefore, the statement follows by Proposition~\ref{prop_double_dual}, from which we know that 
\[\fg\otimes \bU(-1)\cong \Derrho\cong (\Derrho)^{\vee\vee}\cong \Der_\bU^{\chiV}.\]
\end{proof}

\subsection{Stable periodic gradings and Saito matrices}
\label{sec-stable-periodic}

Let us move to the next case of interest, namely cyclic gradings of Lie algebras. We work within the setting of \S\ref{section:theta-representations}.
Let \(\LL\) be a connected semisimple algebraic group and let \(\theta\) be an automorphism of \(\LL\) of finite order \(m\).
Recall that we fixed \(\xi\), a primitive \(m\)-th root of unity. Then we have an \(m\)-cyclic grading of \(\fL\) and a \(G\)-representation \(V\), with:
\[
\fL=\bigoplus_{a\in\ZZ/m\ZZ}\fL_a, \qquad G=(\LL^\theta)^\circ, \qquad V = \fL_1.
\]
We use the bilinear form \(B\) fixed above to identify \(\fL\) with \(\fL^\vee\). Thus, for a polynomial \(F\in\kk[\fL]\), its gradient \(\nabla F:\fL\to\fL\) is defined by setting, for all \(x,y\in\fL\):
\[
B\bigl(\nabla F(x),y\bigr)=d_x F(y).
\]
Consider homogeneous basic invariants of the algebra \(\kk[\fL]^\LL\):
\(
\Phi_1,\ldots,\Phi_r\in\kk[\fL]^\LL,
\)
chosen as eigenvectors for the induced action of \(\theta\) on functions.
For \(1\le i\le r\), write
\[
\deg(\Phi_i)=d_i,
\qquad
a_i\in\ZZ/m\ZZ,
\qquad
\theta\cdot\Phi_i=\xi^{a_i}\Phi_i.
\]
Here our convention for the action on functions is
\(
(\theta\cdot F)(x)=F(\theta^{-1}x).
\)
For \(q\in\ZZ/m\ZZ\), define
\[
I_q:=\{i\in\{1,\ldots,r\}\mid d_i+a_i\equiv q\pmod m\}.
\]

\begin{proposition}
\label{prop-periodic-gradients}
With the notation above:
\begin{enumerate}[label=\roman*)]
\item\label{periodic-restriction-vanishing}
For every \(i\in\{1,\ldots,r\}\), the restriction \(\Phi_i|_V\) vanishes unless \(i\in I_0\).
\item\label{periodic-restriction-generators}
If the grading is stable, the non-zero restrictions \(\Phi_i|_V\), \(i\in I_0\), generate \(\kk[V]^G\).
\item\label{periodic-gradient-degree}
For every \(i\in\{1,\ldots,r\}\),
\[
(\nabla\Phi_i)(V)\subset\fL_{d_i+a_i-1}.
\]
In particular, for \(i\in I_2\),
\[
\eta_i:=(\nabla\Phi_i)|_V:V\to V
\]
is an element of \(\Der_\bU^G\), homogeneous of degree \(d_i-1\).
\item\label{periodic-gradient-basis}
If the grading is stable and \(x\in V\) is stable, then \(x\) is regular semisimple in \(\fL\), the space
\[
\fh=\zeta_{\fL}(x)\cap V
\]
is a Cartan subspace, and the family \((\eta_i(x)\mid i\in I_2)\) forms a basis of \(\fh\). In particular,
\[
|I_2|=\dim_{\kk}(\fh).
\]
\end{enumerate}
\end{proposition}

\begin{proof}
We first prove \ref{periodic-restriction-vanishing}. Fix \(i\in\{1,\ldots,r\}\) and \(x\in V\). Since \(\theta^{-1}x=\xi^{-1}x\), homogeneity and the convention for the action on functions give
\[
\xi^{-d_i}\Phi_i(x)
=\Phi_i(\theta^{-1}x)
=(\theta\cdot\Phi_i)(x)
=\xi^{a_i}\Phi_i(x).
\]
If \(d_i+a_i\not\equiv0\pmod m\), the two scalars are distinct, and therefore \(\Phi_i(x)=0\). Since this holds for every \(x\in V\), the restriction \(\Phi_i|_V\) vanishes. Assertion~\ref{periodic-restriction-generators} is precisely \cite[Proposition~2.29]{laga-romano}.

We next prove \ref{periodic-gradient-degree}. Denote the chosen invariant bilinear form by \(B\). From
\(
\Phi_i\circ\theta=\xi^{-a_i}\Phi_i
\)
and the \(\theta\)-invariance of \(B\), for all \(x,y\in\fL\) we obtain
\[
\begin{aligned}
B\bigl(\nabla\Phi_i(\theta x),\theta y\bigr)
&=d_{\theta x} \Phi_i(\theta y)=d_x (\Phi_i\circ\theta)(y)=\xi^{-a_i}d_x \Phi_i(y)=B\bigl(\xi^{-a_i}\theta\nabla\Phi_i(x),\theta y\bigr).
\end{aligned}
\]
The non-degeneracy of \(B\) therefore gives
\[
\nabla\Phi_i(\theta x)=\xi^{-a_i}\theta\nabla\Phi_i(x).
\]
Now let \(x\in V\). Since \(\theta x=\xi x\) and \(\nabla\Phi_i\) is homogeneous of degree \(d_i-1\), it follows that
\[
\theta\nabla\Phi_i(x)
=\xi^{d_i+a_i-1}\nabla\Phi_i(x).
\]
This proves the stated containment. Moreover, the \(\LL\)-invariance of \(\Phi_i\) and of \(B\) implies that \(\nabla\Phi_i\) is \(\LL\)-equivariant. Hence, when \(i\in I_2\), its restriction \(\eta_i:V\to V\) is \(G\)-equivariant and defines an element of \(\Der_\bU^G\) of degree \(d_i-1\).

Finally, assume that the grading is stable and let \(x\in V\) be stable. By \cite[Proposition~2.18]{laga-romano}, \(x\) is regular semisimple in \(\fL\). Hence
\[
\ft:=\zeta_{\fL}(x)
\]
is a \(\theta\)-stable Cartan subalgebra, and \(\fh=\ft\cap V\) is the associated Cartan subspace.
Kostant’s regularity criterion for the adjoint quotient from \cite[Theorem~9]{konstant:polynomialrings} says that, since \(x\) is regular, the differentials \((d_x \Phi_i\mid 1\le i\le r)\)
are linearly independent. 

Under the identification \(\fL \simeq \fL^\vee\) induced by \(B\), the corresponding gradients are therefore linearly independent. Thus 
\(
(\nabla\Phi_i(x)\mid 1\le i\le r)
\)
form a basis of \(\ft\). Indeed, let us check that \(\nabla\Phi_i(x)\) lies in \(\ft\),   namely, that \([\nabla\Phi_i(x),x]=0\).
Since \(B\) is non-degenerate, it is enough to check
that \(B([\nabla\Phi_i(x),x],y)=0\) for all \(y \in \fL\).
By invariance of \(B\), it suffices to show 
\(B(\nabla\Phi_i(x),[x,y])=0\) for all \(y \in \fL\), which is equivalent to \(d_x\Phi_i([x,y])=0\) for all \(y \in \fL\). But \(\{[x,y]\mid y \in \fL\}\) is the tangent space to the adjoint orbit of \(x\), and \(\Phi_i\) is constant on that orbit. Hence \(d_x\Phi_i([x,y])=0\) for all \(y \in \fL\). 
Now, since \(\dim(\ft)=r\) and \((\nabla\Phi_i(x)\mid 1 \le i \le r)\) are linearly independent by Kostant's theorem on adjoint actions mentioned above, we have the desired basis. Assertion~\ref{periodic-gradient-degree} shows that this basis is compatible with the \(\theta\)-eigenspace decomposition of \(\ft\). Its degree-one part is therefore precisely
\(
(\eta_i(x)\mid i\in I_2).
\)
This is a basis of
\(
\ft_1=\zeta_{\fL}(x)\cap V=\fh.
\)
This proves \ref{periodic-gradient-basis}.
\end{proof}

\begin{proposition} \label{prop-periodic-coexponents} Assume that the grading is stable. Then 
\[ \Coexp(V)=\{d_i -1\mid i\in I_2\}. \] 
Moreover, the restrictions to a Cartan subspace of the vector fields \((\eta_i\mid i\in I_2)\) form a homogeneous \(\bR\)-basis of \(\Der_\bS^W\) and their images under \(\pi_G\) form a basis of \(\Der_\bR(-\log\delta)\). 
\end{proposition}

\begin{proof}
Let \(W_{\fL}\) be the Weyl group of \(\fL\). Its basic invariants
have degrees \(d_1,\ldots,d_r\). Since \(W_{\fL}\) is a real
reflection group, its coexponents are
\(
d_1-1,\ldots,d_r-1.
\)
We use the derivations \(\eta_i \mid i \in \{1,\ldots,r\}\) of Proposition \ref{prop-periodic-gradients}.
Choose a stable element \(x\in V \subset \fL\). By
Proposition~\ref{prop-periodic-gradients}, setting
\(
\ft:=\zeta_{\fL}(x),
\)
the space \(\ft\) is a \(\theta\)-stable Cartan subalgebra of \(\fL\)
and
\[
\fh:=\ft\cap V
\]
is a Cartan subspace. Since \(x\) is regular semisimple in \(\fL\), it lies on no reflecting
hyperplane of \(W_{\fL}\), which means that the triple $(W_{\fL},\theta,\xi)$ is regular in the sense of Lehrer-Springer (we refer to \cite{lehrer-springer:subquotients} and use the terminology \textit{LS-regular} for such a triple). 

By \cite[Corollary~7.4]{reeder-levy-yu-gross}, stability identifies
the Weyl group
\[
W=N_G(\fh)/Z_G(\fh)
\]
with the fixed-point subgroup \(W_{\fL}^{\theta}\) for the induced
action of \(\theta\) on \(W_{\fL}\).

By Proposition \ref{prop-surjadjoint} and the isomorphism in \eqref{DerG}, we have \(
\Der_{\kk[\fL]}^{L} \simeq \Der_{\kk[\ft]}^{W_{\fL}} \simeq  (\kk[\ft] \otimes \ft)^{W_{\fL}}\).
For all \(i \in \{1,\ldots,r\}\), we put \(\Theta_i = \nabla \Phi_i|_\ft\). Then, \((\Theta_i \mid i \in \{1,\ldots,r\})\) is a
homogeneous basis of \(\Der_{\kk[\ft]}^{W_{\fL}}\) over
\(\kk[\ft]^{W_{\fL}}\).
Also, again by the isomorphism \eqref{DerG}, we have \((\bS \otimes \fh)^W \simeq \Der_{\bS}^W\).
Now, by \cite[Theorem D]{lehrer-springer:subquotients}, we obtain the surjectivity of the restriction map 
\[
(\kk[\ft] \otimes \ft)^{W_{\fL}} \to (\bS \otimes \fh)^W
\]
This map sends a \(W_\fL\)-invariant element \(F \in \kk[\ft] \otimes \ft\), seen as \(W_\fL\)-equivariant polynomial map \(\ft \to \ft\) to \(p_1 \circ F|_\fh : \fh \to \fh\), where \(p_1 : \ft \to \fh\) is the projection onto the \(\xi\)-eigenspace \(\fh = \ft \cap V\). Furthermore, the previous maps fit in the following commutative diagram
$$
\xymatrix{
\Der_{\kk[\ft]}^{W_{\fL}} \ar[d] \ar[r] & \Der_\bU^G \ar[d]^{\pi_G}\\
(\kk[\ft] \otimes \ft)^{W_{\fL}} \ar[r] & (\bS \otimes \fh)^W
}
$$
which implies that, in this setting, the restriction of invariant derivations defined by $\pi_G$ is surjective.

Now recall that, by definition,
\(
I_2
=
\{i\in\{1,\ldots,r\}\mid d_i+a_i-1\equiv1\pmod m\}.
\)
For \(i\in I_2\), write
\(\bar\eta_i=\Theta_i|_\fh=\eta_i|_\fh\).
By the surjectivity above and the
eigenspace calculation of Proposition \ref{prop-periodic-gradients}, we have that
\((\bar\eta_i\mid i\in I_2)\) generate \(\Der_\bS^W\)
over \(\bR\).
So, applying the argument of
\cite[Proposition~4.8]{lehrer-springer:subquotients}
to these generators gives the multiset inclusion
\(
\Coexp(W)\subseteq\{d_i-1\mid i\in I_2\}.
\)
However, by Proposition~\ref{prop-periodic-gradients},
part~\ref{periodic-gradient-basis}, we have
\(|I_2|=\dim_\kk(\fh)\). Both multisets therefore have
the same cardinality, so
\[
\Coexp(V)=\Coexp(W)=\{d_i-1\mid i\in I_2\}.
\]
This proves the first statement.

\medskip
Next, choose a stable element \(x\in\fh\). Using again Proposition \ref{prop-periodic-gradients}
part~\ref{periodic-gradient-basis}, we see that the family
\[
(\overline{\eta}_i(x)\mid i\in I_2)
=(\eta_i(x)\mid i\in I_2)
\]
forms a basis of \(\fh\). Hence the \((\overline{\eta}_i \mid i \in I_2)\) are linearly
independent over the fraction field of \(\bR\).

By Theorem~\ref{thm_terao_der0}, the \(\bR\)-module
\(\Der_\bS^W\) has a homogeneous basis
\((\tau_i\mid 1\le i\le n)\) whose degree multiset is \(\Coexp(V)\).
Order this basis and the \(\overline{\eta}_i\) so that their degrees
agree as multisets, and write
\[
(\overline{\eta}_i\mid i\in I_2)
=(\tau_i\mid 1\le i\le n)A
\]
for a matrix \(A\) with entries in \(\bR\). The determinant of \(A\) is
non-zero by the generic independence just proved, while
\[
\deg\bigl(\det(A)\bigr)
=
\sum_{i\in I_2}(d_i-1)
-
\sum_{c\in\Coexp(V)}c
=0.
\]
Since \(\bR = \bigoplus_{t \ge 0} \bR_t\) and \(\bR_0=\kk\), we have
\(\det(A)\in\kk^\times\). Thus \(A\) is invertible over \(\bR\), and
the \(\overline{\eta}_i\) form a homogeneous \(\bR\)-basis of
\(\Der_\bS^W\).
Finally, as established in the proof of Lemma~\ref{lemma-surj},
restriction is compatible with the maps to the quotient:
\[
\pi_G(\eta_i)=\pi_W(\overline{\eta}_i).
\]
The map \(\pi_W\) is an isomorphism by
Theorem~\ref{thm_terao_der0}. So
\((\pi_G(\eta_i)\mid i\in I_2)\) is an \(\bR\)-basis of
\(\Der_\bR(-\log\delta)\).
This proves the second statement.
\end{proof}

\begin{remark}
\label{rmk-inner-congruences}
If the grading is inner, then \(\theta=\operatorname{Int}(g)\) for some \(g\in\LL\). Since every \(\Phi_i\) is \(\LL\)-invariant, the induced action of \(\theta\) on \(\kk[\fL]^\LL\) is trivial. Thus \(a_i=0\) for every \(i\in\{1,\ldots,r\}\), and
\[
I_0=\{i\in\{1,\ldots,r\}\mid d_i\equiv0\pmod m\},
\qquad
I_2=\{i\in\{1,\ldots,r\}\mid d_i\equiv2\pmod m\}.
\]
The first congruence selects the ambient invariants which restrict to basic invariant functions on \(V\), while the second selects the ambient gradients which give invariant derivations. For \(m=2\) these congruences coincide; for larger \(m\) they generally do not. Notice that, when the grading is inner, \cite[Theorem~C(ii)]{lehrer-springer:subquotients} essentially implies the statement above. In \textit{loc. cit.}, the \emph{untwisted} case corresponds precisely to the case when the grading is inner.
\end{remark}

For the next result, we
choose a basis \((e_1,\ldots,e_s)\) of \(\fg\), where
\(s=\dim_{\kk}(\fg)\), and a basis \((v_1,\ldots,v_\ell)\) of \(V\). Let
\((x_1,\ldots,x_\ell)\) be the coordinates dual to the \(v_j\), and let
\((w_1,\ldots,w_\ell)\) be the basis of \(\fL_{-1}\) which is
\(B\)-dual to the \((v_1,\ldots,v_\ell)\).
We consider the coefficient matrix \(M_\psi\) of \(\psi\) and the matrix
\(M_\Phi\) of restricted derivatives of the invariants indexed by
\(I_2\):
\[
M_\psi=
\bigl(\psi(e_a)(x_j)\bigr)_
{1\le a\le s,\ 1\le j\le\ell},
\qquad
M_\Phi=
\left(
\left.\partial_{w_j}\Phi_i\right|_V
\right)_{i\in I_2,\ 1\le j\le\ell}.
\]

\begin{theorem}
\label{thm-periodic-saito}
Assume the grading is stable and that \(G\) is semisimple. Then
\[
\Der_\bU(-\log\Delta)
=
\Derrho\oplus\bigoplus_{i\in I_2}\bU\eta_i
\simeq
\fg\otimes\bU(-1)
\oplus
\bigoplus_{i\in I_2}\bU(1-d_i).
\]
Explicitly, we have a Saito matrix \(M_\Delta\), of the form
\[
M_\Delta=
\begin{pmatrix}
M_\psi\\
M_\Phi
\end{pmatrix},
\]
satisfying
\(\det(M_\Delta)=c\,\Delta\) for some \(c\in\kk^\times\).
\end{theorem}

\begin{proof}
By Proposition~\ref{vinberg-is-polar},
the Vinberg representation \(V\) of \(G\) is visible and polar.
By definition, stability of the grading means that this representation
is stable in the sense of Dadok--Kac and locally free.
In particular, its generic
stabilizer is finite. Since \(G\) is semisimple,
Theorem~\ref{thm_finite_stab} applies and gives the desired decomposition.

We identify explicitly the two summands occurring in that theorem.
Finiteness of the generic stabilizer implies that the infinitesimal
action
\(
\psi:\fg\otimes\bU(-1)\to\Der_\bU
\)
is generically injective, and hence injective. Its image is
\(\Derrho\), so
\(
\Derrho\simeq\fg\otimes\bU(-1).
\)
For the remaining summand, Proposition~\ref{prop-periodic-coexponents}
shows that
\(
\bigl(\pi_G(\eta_i)\mid i\in I_2\bigr)
\)
is an \(\bR\)-basis of \(\Der_\bR(-\log\delta)\).
Therefore, in the proof of Theorem~\ref{thm_direct_sum_dec},
we may choose the
\(\bR\)-linear section defined on this basis by
\[
\sigma\bigl(\pi_G(\eta_i)\bigr)=\eta_i,
\qquad
\text{for every }i\in I_2.
\]
After extension of scalars from \(\bR\) to \(\bU\), the image of this
section is exactly \(\bigoplus_{i\in I_2}\bU\eta_i\). The split exact
sequence used in that theorem now gives the explicit direct-sum
decomposition
\[
\Der_\bU(-\log\Delta)
=
\Derrho\oplus\bigoplus_{i\in I_2}\bU\eta_i.
\]
Since \(\deg(\eta_i)=d_i-1\) for every \(i\in I_2\), this also gives the displayed graded
isomorphism.

It remains to identify the coefficient matrix of this basis. By
definition, for \(1\le a\le s\), the \(a\)-th row of \(M_\psi\) consists of the coefficients of
the infinitesimal-action vector field \(\psi(e_a)\) in the basis of
coordinate vector fields.
Consequently the rows of \(M_\psi\) are the coefficient rows of the above
basis of \(\Derrho\). For \(i\in I_2\) and \(1\le j\le\ell\), let \(c_{ij}\in\bU\) be the
coefficient functions determined by
\[
\eta_i(x)=\sum_{j=1}^{\ell}c_{ij}(x)v_j,
\qquad x\in V.
\]
Because \(w_1,\ldots,w_\ell\) is \(B\)-dual to
\(v_1,\ldots,v_\ell\), the definition of the gradient gives
\[
c_{ij}(x)
=
B\bigl(\nabla\Phi_i(x),w_j\bigr)
=
d_x\Phi_i(w_j)
=
\partial_{w_j}\Phi_i(x).
\]
Restricting to \(x\in V\) shows that the coefficient matrix of the
basis \((\eta_i\mid i\in I_2)\) of the remaining summand is precisely
\(M_\Phi\).

If \(n=\dim_{\kk}(\fh)\), then
\(|I_2|=n\) by
Proposition~\ref{prop-periodic-gradients},
part \ref{periodic-gradient-basis}. Moreover, the generic stabilizer
is finite, so a general orbit has dimension \(\dim_{\kk}(\fg)\); polarity
gives
\[
\ell=\dim_{\kk}(\fg)+n.
\]
Thus \(M_\psi\) has \(\dim_{\kk}(\fg)\) rows, \(M_\Phi\) has \(n\) rows, and the
stacked matrix \(M_\Delta\) is square of size \(\ell\). The preceding
direct-sum decomposition shows that its rows form a
\(\bU\)-basis of \(\Der_\bU(-\log\Delta)\).

By Lemma~\ref{lem-polar-discriminant-reduced}, \(\Delta\) is reduced.
Saito's criterion applied to the above basis therefore yields
\[
\det(M_\Delta)=u\,\Delta
\]
for a unit \(u\in\bU^\times\). Since \(\bU\) is a polynomial ring over
\(\kk\), its units are the elements of \(\kk^\times\); writing
\(u=c\) proves the last assertion.
\end{proof}

\begin{remark}
The gradient construction in Proposition~\ref{prop-periodic-gradients} does not require the degree-zero group to be semisimple. Thus it also applies to stable periodic gradings for which \(G\) has a torus factor, including Coxeter gradings. The semisimplicity assumption in Theorem~\ref{thm-periodic-saito} is used both in Proposition~\ref{prop_double_dual} and in Lemma~\ref{lem-polar-discriminant-reduced}. Weakening it would require a separate analysis of vertical derivations and of the reducedness of the polar discriminant for reductive degree-zero groups; see Remark~\ref{rmk-order6-e8-reducedness}.
\end{remark}

\subsection{Exceptional examples}

Discriminant divisors obtained from cyclic gradings of exceptional simple Lie algebras of type \(\fe_6\),  \(\fe_7\) and \(\fe_8\) are particularly interesting and well-studied.
Stable gradings of exceptional simple Lie algebras are classified in
\cite[\S8.1]{reeder-levy-yu-gross}. We record five gradings, four of which are defined by an inner automorphism and a fifth by an outer one, and apply our results to give an explicit Saito matrix in each case.

Some further examples of outer gradings arising from exceptional automorphisms provide 
a 3-grading of \(\fso_8\) inherited from triality, and a 4-grading of \(\fe_6\). The resulting representations are associated with elliptic curves 
of degree 3 and 4 as they give, respectively, \(\fsl_3\) acting on \(S^3 V_3\) and 
\(\fsl_2 \oplus \fsl_4\) acting on \(A_2 \otimes S^2 B_4\), where \(V_3\), \(A_2\) and \(B_4\) are vector spaces of dimension \(3\), \(2\) and \(4\).
These toy cases are related, respectively, to work of Bhargava-Ho, \cite{bhargava-ho:genus-one} and Buchweitz-Pike,  \cite{buchweitz-pike}. 
In the first case, we get a \(10\times 10\) Saito matrix for the classical degree-\(12\) discriminant of plane cubics, where, besides the Euler derivation, the remaining invariant basis derivation arises from the Hessian.
In the second case, one gets a free divisor of degree \(24\)
as determinant of a Saito matrix of size \(20 \times 20\).

\subsubsection{The stable \(2\)-grading of \(\fe_7\)}
\label{2-cyclic of e7}
We consider first
\[
\fL=\fe_7,
\qquad
\fg\simeq\fsl(V_8),
\qquad
V\simeq\wedge^4 V_8,
\]
so \(\dim_{\kk}(\fg)=63\) and \(\dim_{\kk}(V)=70\). The Weyl group \(W\) is \(W(E_7)\), equivalently \(G_{36}\) in the
Shephard--Todd classification. The ambient basic invariant degrees are
\[
2,6,8,10,12,14,18.
\]
Since \(m=2\), they all lie in both \(I_0\) and \(I_2\), and
\[
\Coexp(V)=\{1,5,7,9,11,13,17\}.
\]
Therefore we get
\[
\Der_\bU(-\log\Delta)
\simeq
\fsl_8\otimes\bU(-1)
\oplus
\bigoplus_{c\in\{1,5,7,9,11,13,17\}}\bU(-c).
\]
The Saito determinant has degree \(63+(1+5+7+9+11+13+17)=126\).

\subsubsection{The stable \(2\)-grading of \(\fe_8\)}
\label{2-cyclic of e8}
Next, we look at the \(2\)-cyclic grading of \(\fe_8\), which gives
\[
\fL=\fe_8,
\qquad
\fg\simeq\fso(V_{16}),
\qquad
V\simeq\Sigma_+,
\]
where \(\Sigma_+\) is a half-spin representation. Thus \(\dim_{\kk}(\fg)=120\) and \(\dim_{\kk}(V)=128\). The Weyl group \(W\) is
\(W(E_8)\), equivalently the Shephard--Todd group \(G_{37}\), and
\[
\Coexp(V)=\{1,7,11,13,17,19,23,29\}.
\]
Hence
\[
\Der_\bU(-\log\Delta)
\simeq
\fso_{16}\otimes\bU(-1)
\oplus
\bigoplus_{c\in\{1,7,11,13,17,19,23,29\}}\bU(-c).
\]
The Saito determinant has degree \(120+120=240\) and gives the degree-\(240\) discriminant on \(\Sigma_+\).

\subsubsection{The stable \(3\)-grading of \(\fe_8\)}
\label{3-cyclic of e8}
The next case of interest is the \(3\)-cyclic grading of \(\fe_8\). It gives 
\[
\fL=\fe_8,
\qquad
\fg\simeq\fsl(V_9),
\qquad
V\simeq\wedge^3V_9,
\]
so \(\dim_{\kk}(\fg)=80\) and \(\dim_{\kk}(V)=84\). The ambient \(E_8\) degrees split as
\[
\{d_i\mid i\in I_0\}=\{12,18,24,30\},
\qquad
\{d_i\mid i\in I_2\}=\{2,8,14,20\}.
\]
The invariants indexed by \(I_0\) restrict to basic invariant functions
on \(V\), while the gradients indexed by \(I_2\) give invariant
derivations of degrees
\(
1,7,13,19.
\)
The Weyl group \(W\) is the Shephard--Todd group \(G_{32}\), and these
are its coexponents. Consequently we get
\[
\Der_\bU(-\log\Delta)
\simeq
\fsl_9\otimes\bU(-1)
\oplus \bigoplus_{c \in \{1,7,13,19\}} \bU(-c)
.
\]
The infinitesimal-action block has size \(80\times84\), the gradient block has size \(4\times84\), and the Saito determinant has degree \(120\).

\subsubsection{The stable \(5\)-grading of \(\fe_8\)}
\label{5-cyclic of e8}
This is the Vinberg representation used in the invariant-theoretic description of \(5\)-Selmer elements; see \cite[\S5]{laga-romano}. Up to a finite central quotient,
\[
G\simeq\SL(A_5)\times\SL(B_5),
\qquad
V\simeq\wedge^2A_5\otimes B_5,
\]
where \(A_5\) and \(B_5\) are \(5\)-dimensional vector spaces, so \(\dim_{\kk}(\fg)=48\) and \(\dim_{\kk}(V)=50\). For the ambient \(E_8\) degrees
\[
2,8,12,14,18,20,24,30,
\]
one has
\(
\{d_i\mid i\in I_0\}=\{20,30\},
\) and \(
\{d_i\mid i\in I_2\}=\{2,12\}.
\)
Thus the ambient invariants of degrees \(20\) and \(30\) restrict to the
two basic invariant functions, while the gradients of the ambient
invariants of degrees \(2\) and \(12\) give invariant derivations of
degrees \(1\) and \(11\). The Weyl group \(W\) is the Shephard--Todd group \(G_{16}\), and these
are its coexponents. Therefore
\[
\Der_\bU(-\log\Delta)
\simeq
(\fsl_5 \oplus \fsl_5)\otimes\bU(-1)\oplus\bU(-1)\oplus\bU(-11),
\]
and the Saito determinant has degree \(48+1+11=60\).

\subsubsection{The stable outer \(2\)-grading of \(\fe_6\)}
\label{outer of e6}
The affine Dynkin diagram associated to the outer involution of $E_6$ is the one of $F_4$ to which a node is added next to the fourth (short) simple root (refer to \cite{Vinberg}). Among the gradings one can construct from outer automorphisms of $E_6$, we choose the one appearing in \cite{benedetti-manivel:discriminants}, so that we can apply Proposition \ref{prop_thetacorank_stable}. Up to a finite central quotient,
\[
G\simeq\Sp(V_8),
\qquad
V\simeq\wedge^{\langle 4\rangle }V_8,
\]
where $\wedge^{\langle 4\rangle }V_8=\wedge^4V_8/(\Omega\wedge\wedge^2V_8)$, with $\Omega\in\wedge^2V_8$ the invariant bivector corresponding to the symplectic form under the identification $V_8\simeq V_8^\vee$. We have \(\dim_{\kk}(\fg)=36\) and \(\dim_{\kk}(V)=42\). The ambient \(E_6\) degrees are
\[
2,5,6,8,9,12.
\]
For this stable outer involution, the Weyl group \(W\) is
\(W(E_6)\), and
\[
\Coexp(V)=\{1,4,5,7,8,11\}.
\]
Therefore we get
\[
\Der_\bU(-\log\Delta)
\simeq
\fsp_8\otimes\bU(-1)\oplus\bigoplus_{c \in \{1,4,5,7,8,11\}}\bU(-c).
\]
The Saito determinant has degree \(36+1+4+5+7+8+11=72\). Notice that the degree-$72$ discriminant in this case is defined by the restriction of the $E_6$-discriminant equation to $\wedge^{\langle 4\rangle }V_8$ and, by \cite[Theorem 5.2]{benedetti-manivel:discriminants}, it identifies with the dual projective variety of the Lagrangian Grassmannian $\mathrm{LG}(4,8)$.

\bibliographystyle{amsalpha.v2}
\bibliography{discriminants}

\end{document}